\documentclass[11pt, eqno]{article}
\usepackage{amssymb,amsmath,latexsym}
\usepackage{epsfig}
\usepackage{color}

\newtheorem{thm}{Theorem}[section]
\newtheorem{cor}[thm]{Corollary}
\newtheorem{lemma}[thm]{Lemma}
\newtheorem{prop}[thm]{Proposition}
\newtheorem{defn}[thm]{Definition}

\newtheorem{remark}[thm]{Remark}
\numberwithin{equation}{section}

\def\pf{{\medskip\noindent {\bf Proof. }}}
\def\qed{{\hfill $\Box$ \bigskip}}

\def\sA {{\cal A}}  \def\sC {{\cal C}}
 \def\sE {{\cal E}} \def\sF {{\cal F}}

\def\sM {{\cal M}}

 \def\bH {{\mathbb H}}

  \def\bR {{\mathbb R}}

\def\R {{\mathbb R}}

\def\wt{\widetilde}
\def\wh{\widehat}

\def\E{{\mathbb E}}
\def\P{{\mathbb P}}

\def\bea{\begin{align*}}
\def\eea{\end{align*}}
\def\bee{\begin{equation}}
\def\eee{\end{equation}}

\def\eps{\varepsilon}

\def\wh{\widehat}

\def\1{{\bf 1}}

\def\loc{{\rm loc}}

\makeatletter
\@addtoreset{equation}{section}

\makeatother

\begin{document}
\bibliographystyle{plain}

\title{\Large \bf
Boundary Harnack principle for
    $\Delta + \Delta^{\alpha/2}$}

\author{{\bf Zhen-Qing Chen}\thanks{Research partially supported
 by NSF Grants DMS-0600206 and DMS-0906743.},
\quad {\bf Panki Kim}\thanks{Research
  supported by National Research Foundation of Korea Grant funded by
the Korean Government (2009-0087117)}, \quad {\bf Renming Song} \quad and \quad
{\bf Zoran Vondra\v{c}ek}\thanks{Research partially supported by the
MZOS grant 037-0372790-2801 of the Republic of Croatia.}}

\date{(November 8, 2009)}

\maketitle

\begin{abstract}
For $d\geq 1$ and $\alpha \in (0, 2)$, consider the family of pseudo
differential operators $\{\Delta+ b \Delta^{\alpha/2}; b\in [0,
1]\}$ on $\R^d$ that evolves continuously from $\Delta$ to $\Delta +
\Delta^{\alpha/2}$. In this paper, we establish a uniform boundary
Harnack principle (BHP) with explicit boundary decay rate for
nonnegative functions which are harmonic with respect to $\Delta +b
\Delta^{\alpha/2}$ (or equivalently, the sum of a Brownian motion
and an independent symmetric $\alpha$-stable process with constant
multiple $b^{1/\alpha}$) in $C^{1, 1}$ open sets. Here a ``uniform"
BHP  means that the comparing constant in the BHP is independent of
$b\in [0, 1]$. Along the way, a uniform Carleson type estimate is
established for nonnegative functions which are harmonic with
respect to $\Delta + b \Delta^{\alpha/2}$  in Lipschitz open sets.
Our method employs a combination of probabilistic and analytic
techniques.
\end{abstract}



\bigskip
\section{Introduction}

Discontinuous Markov processes have been receiving intensive study
recently due to their importance both in theory and in applications.
Many physical and economic systems could be and in fact have been
successfully modeled by discontinuous Markov processes (or jump
diffusions as some authors call them); see for example, \cite{JW,
KSZ, OS} and the references therein. The infinitesimal generator of
a discontinuous Markov process in $\R^d$ is no longer a differential
operator but rather a non-local (or integro-differential) operator.
For instance, the infinitesimal generator of a rotationally
symmetric $\alpha$-stable process in $\R^d$ with $\alpha \in (0, 2)$
is a fractional Laplacian operator $c\, \Delta^{\alpha /2}:=- c\,
(-\Delta)^{\alpha /2}$.

Discontinuous Markov processes include the very important L\'evy
processes as special cases and they are of intrinsic importance in
probability theory. Integro-differential operators are very
important in the theory of partial differential equations. Most of
the recent study concentrates on discontinuous Markov processes,
like the rotationally symmetric $\alpha$-stable processes, that do
not have a diffusion component. For a summary of some of these
recent results from the probability literature, one can see
\cite{C0, BBKRSV} and the references therein. We refer the readers
to \cite{CSS, CaS, CV} for a sample of
 recent progresses in the PDE literature.

However, in many situations, like in finance and control theory, one
needs Markov processes that have both a diffusion component and a
jump component, see for instance,
\cite{JKC, MK2, OS}.
The fact that such a process $X$ has both
diffusion and jump components is the source of many difficulties in
investigating the potential theory of the process $X$. The main
difficulty in studying $X$ stems from the fact that it runs on two
different scales: on the small scale the diffusion part dominates,
while on the large scale the jumps take over. Another difficulty is
encountered when looking at the exit of $X$ from an open set: for
diffusions, the exit is through the boundary, while for the pure
jump processes, typically the exit happens by jumping out from the
open set. For the process $X$, both cases will occur which makes the
process $X$ much more difficult to study.

Despite these difficulties, in the last few years significant
progress has been made in understanding the potential theory of such
processes. Green function estimates (for the whole space) and the
Harnack inequality for a class of processes with both continuous and
jump components were established in \cite{RSV} and \cite{SV05}. The
parabolic Harnack inequality and heat kernel estimates were studied
in \cite{SV07} for L\'evy processes on $\bR^d$ that are the
independent sum of Brownian motion and symmetric stable process, and
in \cite{CK08} for much more general symmetric diffusions with
jumps. Moreover,  a priori H\"older estimate is established in
\cite{CK08} for bounded parabolic functions. For earlier results on
second order integro-differential operators, one can see \cite{GM}
and the references therein.

The boundary Harnack principle (BHP) is a result about the ratio of
positive harmonic functions. We say that the BHP holds for an open
set $D\subset \R^d$ if there exist positive constants $R_0$ and $C$
depending on $D$ with the property that for any $Q\in
\partial D$, $r\in (0, R_0]$, and any positive harmonic functions
$u$ and $v$ in $D\cap B(Q, r)$ that vanish continuously on $\partial
D\cap B(Q, r)$, we have
\begin{equation}\label{e:bhp_b}
\frac{u(x)}{u(y)}\le C\frac{v(x)}{v(y)} \qquad \mbox{ for all } x,
y\in D\cap B(Q, r/2).
\end{equation}

The BHP for Brownian motion (or, equivalently, for the Laplacian) is
a fundamental result in analysis and PDE.  It was independently
established for Lipschitz domains in the late 1970's by Ancona,
Dahlberg and Wu (\cite{An1,Da,W}). Later, Bass and Burdzy developed
a probabilistic method in \cite{BB2} to prove the boundary Harnack
principle and extended the boundary Harnack principle to more
general domains (see also \cite{BB1}). When $D$ is a bounded
$C^{1,1}$ domain, \eqref{e:bhp_b} can be strengthened to the
following version that gives the explicit boundary decay rate of
non-negative harmonic functions that vanish on the boundary:
\begin{equation}\label{e:0.3}
\frac{u(x)}{u(y)}\le C\frac{\delta_D(x)}{\delta_D(y)} \qquad \mbox{ for all } x,
y\in D\cap B(Q, r/2),
\end{equation}
where $\delta_D(x)$ is the Euclidean distance between $x$ and $D^c$.
The BHP plays a vital role in the study of potential theory of
Brownian motion and Dirichlet Laplacian in domains. For example, BHP
can be used to show that Martin boundary can be identified with the
Euclidean boundary for a large class of domains and to study the
non-tangential limit of non-negative harmonic functions near the
boundary (see \cite{AG} for an analytic approach and \cite{Ba} for a
probabilistic approach). In fact, BHP has also be established for a
large class of diffusions (or, equivalently, for second order
elliptic equations), see \cite{CFMS, FGMS}.

The study of BHP for discontinuous Markov processes and
integro-differential operators is quite recent. It was first
established for bounded Lipschitz domains in \cite{Bo} and then
extended to more general open sets in \cite{SW99}.
 Subsequently Bogdan-Stos-Sztonyk \cite{BSS} and Sztonyk \cite{Sz2}
extended the boundary Harnack principle to symmetric (but not
necessarily rotationally invariant) stable processes. Recently, the
BHP has been extended in \cite{KSV} to a large class of pure jump
L\'evy processes that can be obtained from Brownian motion through
subordination. Very recently, the boundary Harnack principle for
some one-dimensional L\'evy processes with both continuous and jump
components was studied in \cite{KSV09}. However BHP for processes on
$\bR^d$ in dimension two and higher that have both diffusion and
jump components have been  completely  open until now. Note that
 the fact that a pure jump process may (and typically does) exit an open
set by jumping out of it stipulates that, in the boundary Harnack
principle for such processes, the nonnegative harmonic functions
vanish continuously on $D^c \cap B(Q,r)$.

The principal goal of this paper is to establish the boundary
Harnack principle for nonnegative functions which are harmonic with
respect to the independent sum of a Brownian motion and a symmetric
stable process in $C^{1, 1}$ open sets in $\R^d$ for every $d \ge
1$. The process $X$ studied in this paper, although quite specific,
serves as a test case for more general processes with both
continuous and jump parts. The study of this test case will
hopefully shed new light on the understanding of the boundary
behavior of nonnegative harmonic functions of general Markov
processes.

Intuitively, the independent sum $X$ of a Brownian motion and a
symmetric stable process can be thought roughly as some sort of
``perturbation'' of Brownian motion. Thus some people might expect
the BHP for $X$ could be established through some general
perturbation technique. However, this  kind of approach may not
always work. In \cite{KS, KS2}, the potential theory of truncated
symmetric stable processes including BHP was studied. One of the
main results in \cite{KS} is that the BHP is valid for the positive
harmonic functions of this process in bounded convex domains. A very
interesting fact is, even though truncated symmetric stable
processes can be considered as a perturbation  of rotationally
symmetric stable processes (see \cite{GR, KS2}), unlike symmetric
stable processes, the BHP for truncated symmetric stable processes
 fails in
non-convex domains (see the last section of
\cite{KS} for a counterexample). This indicates that general
perturbation method may not be suitable for establishing the BHP.

Let us now describe the main result of this paper more precisely
 and at the same time fix the notations. A (rotationally) symmetric
$\alpha$-stable process $Y=(Y_t, t\geq 0, \P_x, x\in \bR^d)$ in
$\bR^d$ is a L\'evy process such that
$$
\E_x \left[ e^{i\xi\cdot(Y_t-Y_0)} \right]\,=\,e^{-t|\xi|^{\alpha}}
\qquad \hbox{for every } x\in \bR^d \hbox{ and }  \xi\in \bR^d.
$$
The infinitesimal generator of a symmetric $\alpha$-stable process
$Y$ in $\bR^d$ is the fractional Laplacian $\Delta^{\alpha /2}$,
which is a prototype of nonlocal operators. The fractional Laplacian
can be written in the form
\begin{equation}\label{e:1.1}
\Delta^{\alpha /2} u(x)\, =\, \lim_{\eps \downarrow 0}\int_{\{y\in
\bR^d: \, |y-x|>\eps\}} (u(y)-u(x)) \frac{\sA (d, \alpha)
}{|x-y|^{d+\alpha}}\, dy
\end{equation}
where $ {\cal A}(d, \alpha):= \alpha2^{\alpha-1}\pi^{-d/2}
\Gamma(\frac{d+\alpha}2) \Gamma(1-\frac{\alpha}2)^{-1}$. Here
$\Gamma$ is the Gamma function defined by $\Gamma(\lambda):=
\int^{\infty}_0 t^{\lambda-1} e^{-t}dt$ for every $\lambda > 0$.

Suppose $X^0$ is a Brownian motion in $\R^d$ with generator
$\Delta=\sum_{i=1}^d \frac{\partial^2}{\partial x_i^2}$, and $Y$ is
a symmetric $\alpha$-stable process in $\R^d$. Both $X^0$ and $Y$
satisfy a self-similarity, which will be used several times in this
paper. That is, for every $\lambda >0$,
$\{\lambda^{-1/2}(X^0_{\lambda t} -X^0_0), t\geq 0\}$ and
$\{\lambda^{-1/\alpha}(Y_{\lambda t} -Y_0), t\geq 0\}$ have the same
distributions as that of $\{X^0_t-X^0_0, t\geq 0\}$ and $\{Y_t-Y_0,
t\geq 0\}$, respectively. Assume that $X^0$ and $Y$ are independent.
For any $a>0$, we define $X^a$ by $X_t^a:=X^0_t+ a Y_t$. We will
call the process $X^a$ the independent sum of the Brownian motion
$X^0$ and the symmetric $\alpha$-stable process  $Y$ with weight
$a>0$. The infinitesimal generator of $X^a$ is $\Delta+ a^\alpha
\Delta^{\alpha/2}$. For every open subset $D\subset \bR^d$, we
denote by  $X^{a,D}$ the subprocess of $X^a$ killed upon leaving
$D$. The infinitesimal generator of $X^{a,D}$ is $(\Delta+ a^\alpha
\Delta^{\alpha/2})|_D$. It is known (see \cite{SV07}) that $X^{a,D}$
has a continuous transition density $p^a_D(t, x, y)$ with respect to
the Lebesgue measure. We will use $p^a(t, x, y)$ to denote the
transition density of $X^a$ (or equivalently, the heat kernel of
$\Delta + a^\alpha \Delta^{\alpha/2}$). The quadratic form $(\sE,
\sF)$ associated with the generator $\Delta + a^\alpha
\Delta^{\alpha/2}$ of $X^a$ is given by
$$
\sF=W^{1,2}(\R^d):=\left\{u\in L^2(\R^d; dx): \,
 \frac{\partial u}{\partial x_i} \in L^2(\R^d; dx)
 \, \hbox{ for every } \,
  1\le i \le d \right\}
$$
and for  $u, v\in \sF$,
$$
\sE(u, v) = \int_{\R^d} \nabla u(x) \cdot \nabla v(x) \, dx +
\frac{1}{2} \int_{\R^d\times \R^d} (u(x)-u(y))(v(x)-v(y)) \frac{
\sA(d, \alpha)\, a^\alpha }{|x-y|^{d+\alpha}} dxdy .
$$
In probability theory, the quadratic form $(\sE, W^{1,2}(\R^d))$ is
called the Dirichlet form of $X^a$. A statement is said to hold
quasi-everywhere (q.e.~in abbreviation) if there is a set $N$ having
zero capacity with respect to $(\sE_1, W^{1,2}(\R^d))$ such that the
statement holds everywhere outside $N$. Here $\sE_1(u, u):=\sE(u,
u)+\int_{\R^d} u(x)^2 dx$. The function $J^a (x, y):=a^\alpha \sA(d,
\alpha)|x-y|^{-(d+\alpha)}$ is the L\'evy intensity of $X^a$. It
determines the L\'evy system for $X^a$, which describes the jumps of
the process $X^a$: for any non-negative measurable function $f$ on
$\bR_+ \times \bR^d\times \bR^d$, $x\in \bR^d$ and stopping time $T$
(with respect to the filtration of $X^a$),
\begin{equation}\label{e:levy}
\E_x \left[\sum_{s\le T} f(s,X^a_{s-}, X^a_s) \right]= \E_x \left[
\int_0^T \left( \int_{\bR^d} f(s,X^a_s, y) J^a(X^a_s,y) dy \right)
ds \right].
\end{equation}
(see, for example, \cite[Proof of Lemma 4.7]{CK1} and \cite[Appendix A]{CK2}.)

The purpose of this paper is to establish
the scale invariant version of the boundary Harnack principle in Theorem
\ref{t:main}. To state this theorem,  we first recall that an open
set $D$ in $\bR^d$ (when $d\ge 2$) is said to be $C^{1,1}$ if there
exist a localization radius $R>0$ and a constant $\Lambda
>0$ such that for every $Q\in \partial D$, there exist a
$C^{1,1}$-function $\phi=\phi_Q: \bR^{d-1}\to \bR$ satisfying $\phi
(0)= \nabla\phi (0)=0$, $\| \nabla \phi  \|_\infty \leq \Lambda$, $|
\nabla \phi (x)-\nabla \phi (y)| \leq \Lambda |x-y|$, and an
orthonormal coordinate system $CS_Q$: $y=(y_1, \cdots, y_{d-1},
y_d)=:(\wt y, \, y_d)$ with its origin at $Q$ such that
$$
B(Q, R)\cap D=\{ y=(\wt y, y_d)\in B(0, R) \mbox{ in } CS_Q: y_d >
\phi (\wt y) \}.
$$
The pair $(R, \Lambda)$ is called the characteristics of the
$C^{1,1}$ open set $D$. By a $C^{1,1}$ open set in $\bR$ we mean an
open set which can be written as the union of disjoint intervals so
that the minimum of the lengths of all these intervals is positive
and the minimum of the distances between these intervals is
positive. Note that a $C^{1,1}$ open set can be unbounded and
disconnected.

For any $x\in D$, let $\delta_D(x)$ denote the
distance between $x$ and $\partial D$. It is well known that any
$C^{1, 1}$ open set $D$ satisfies the uniform interior ball
condition: there exists $\wt R\le R$ such that for every $x\in D$ with
$\delta_{ D}(x)< \wt R$, there is $Q_x \in
\partial D$ so that $|x-Q_x|=\delta_{D}(x)$ and that
$B(\wt{x},\wt R)\subset D$ for $\wt{x}=Q_x+\wt R(x-Q_x)/|x-Q_x|$.
Without loss of generality, throughout this paper, we assume that
the characteristics $(R, \Lambda)$ of a $C^{1, 1}$ open set
satisfies $R=\wt R\le 1$ and $\Lambda\ge 1$.

For any open set $D\subset \bR^d$, $\tau^a_D:=\inf\{t>0: \,
X^a_t\notin D\}$ denotes the first exit time from $D$ by $X^a$.

\begin{defn}\label{D:1.1}  \rm A real-valued function $u$ defined on
$\R^d$ is said to be harmonic in $D\subset \R^d$ with respect to
$X^a$ if for every open set $B$ whose closure is a compact subset of
$D$,
\begin{equation}\label{e:har}
\E_x \left[ \big| u(X^a_{\tau^a_{B}})\big| \right]<\infty \quad
\hbox{and} \quad u(x)= \E_x \left[ u(X^a_{\tau^a_{B}})\right] \qquad
\hbox{for q.e. } x\in B.
\end{equation}
\end{defn}

Note that by using the L\'evy system of $X^a$, we have
\begin{eqnarray*}
\E_x \left[ \big| u(X^a_{\tau^a_{B}})\big| \right] &\geq& \E_x \left[ \big|
u(X^a_{\tau^a_{B}})\big|; \, X^a_{\tau_B}\in \R^d \setminus \overline B \, \right] \\
&=& \E_x \left[ \int_0^{\tau_B} \left(\int_{\R^d \setminus \overline
B} \, |u(y)| \, \frac{ \sA(d, \alpha) \, a^\alpha
}{|X_s^a-y|^{d+\alpha}}\, dy \right) ds\right].
\end{eqnarray*}
Hence if $u$ is a harmonic function in $D$ with respect to $X^a$,
then $u(y) (1\wedge |y|^{-(d+\alpha)})$ is integrable on $B^c$ for
any relatively compact open subset $B$ with $\overline B\subset D$.
It follows from Theorems 1.2 and 1.3 of \cite{CK08} that all
harmonic functions in $D$ with respect to  $X^a$ are continuous on
$D$, since every harmonic function in $D$ with respect to $X^a$ can
be approximated locally uniformly in $D$ by functions that are
bounded on $\R^d$ and harmonic with respect to $X^a$ in relatively
compact open subsets of $D$. Therefore, for any harmonic function
$u$ in $D$, \eqref{e:har} holds for {\it every} point $x\in D$. The
above also implies that any harmonic function $u$ in $D$ with
respect to $X^a$ is locally bounded in $D$ with $\int_{\R^d}
|u(y)|(1\wedge |y|^{-(d+\alpha)})dy <\infty$.   A function $u$ is
said to be in $W^{1,2}_\loc (D)$ if for every relatively compact
subset $B$ with $\overline B\subset D$, there is a function $f\in
W^{1,2}(\R^d)$ such that $u=f$ a.e. on $B$. The following analytic
characterization of a function $u$ being harmonic in $D$ with
respect to $X^a$ follows immediately from Example 2.14 in \cite{C}.

\begin{prop}\label{L:harmonic}
Let $D$ be an open subset of $\R^d$. Then the following are
equivalent.
\begin{description}
\item{\rm (i)} $u$  is harmonic in $D$ with respect to $X^a$;
\item{\rm (ii)} $u$ is locally bounded in $D$, $\int_{\R^d} |u(y)|
(1\wedge |y|^{-(d+\alpha)})dy<\infty$, $u\in W^{1,2}_\loc (D)$ and
$(\Delta + a^\alpha \Delta^{\alpha/2})u=0$ in $D$ in the
distributional sense: for every $ \phi\in C^\infty_c(D) $
$$
\int_{\R^d} \nabla u(x) \cdot \nabla \phi (x) \, dx+ \frac12
\int_{\R^d\times \R^d} (u(x)-u(y))(\phi(x)-\phi(y)) \frac{\sA(d,
\alpha) \, a^\alpha}{|x-y|^{d+\alpha}} \, dxdy =0 .
$$
\end{description}
\end{prop}

The following uniform Harnack principle will be used to prove the
main result of this paper.

\begin{prop}[Harnack principle]\label{uhp}
Suppose that $M>0$. There exists a constant $C_0=C_0(\alpha, M)>0$
such that for any $r\in (0, 1]$, $a\in [0, M]$, $x_0\in \R^d$ and
any function $u$ which is nonnegative in $\R^d$ and harmonic in
$B(x_0, r)$ with respect to $X^a$ we have
$$
u(x)\le C_0u(y) \qquad \mbox{ for all } x, y\in B(x_0, r/2).
$$
\end{prop}

\medskip

Let $Q\in \partial D$. We will say that a function $u:\bR^d\to \R$
vanishes continuously on $ D^c \cap B(Q, r)$ if $u=0$ on $ D^c \cap
B(Q, r)$ and $u$ is continuous at every point of $\partial D\cap
B(Q,r)$. The following is the main result of this paper.

\begin{thm}\label{t:main}
Suppose that $M>0$. For any $C^{1, 1}$ open set $D$ in $\bR^d$ with
the characteristics $(R, \Lambda)$, there exists a positive constant
$C=C(\alpha, d, \Lambda, R, M)$ such that for $a \in [0, M]$, $r \in
(0, R]$, $Q\in \partial D$ and any nonnegative function $u$ in
$\R^d$ that is harmonic in $D \cap B(Q, r)$ with respect to $X^{a}$
and vanishes continuously on $ D^c \cap B(Q, r)$, we have
\begin{equation}\label{e:bhp_m}
\frac{u(x)}{u(y)}\,\le
C\,\frac{\delta_D(x)}{\delta_D(y)}  \qquad
\hbox{for every } x, y\in  D \cap B(Q, r/2).
\end{equation}
\end{thm}

When $a$ changes from $0$ to $M$, $\Delta +a^\alpha
\Delta^{\alpha/2}$ changes continuously from $\Delta$ to $\Delta
+M^\alpha \, \Delta^{\alpha/2}$. So Theorem \ref{t:main} says that
the BHP holds uniformly for the family $\{ \Delta +a^\alpha
\Delta^{\alpha/2}, \, a\in [0, M]\}$ of pseudo differential
operators in the sense that the constant $C$ in \eqref{e:bhp_m} can
be chosen to be independent of $a\in [0, M]$. Note that $a=0$
corresponds to the classical case of the boundary Harnack principle
for the Laplacian. We will therefore in the rest of the paper assume
that $a\in (0,M]$.

As far as we know, this is the first time that a BHP has been
established for non-local integro-differential operators that have
second order differential operator components in dimension two and
higher. Unlike \eqref{e:bhp_b} and the paragraph following it, in
this paper we are concerned with the above BHP for $C^{1,1}$ open
sets only. The main focus and goal of this paper  is to get the
explicit decay rate of harmonic functions near the boundary of $D$
as in \eqref{e:bhp_m} and to show that the BHP is {\it uniform} in
$a\in [0, M]$. We emphasize that \eqref{e:bhp_m} is not true in
Lipschitz domains even in the classical case of BHP for the
Laplacian. However, a uniform Carleson type estimate is shown to
hold for Lipschitz open sets in Theorem \ref{carleson}. The BHP of
above type is very useful in studying other fine properties of the
process. For example, we will use it to derive sharp Green function
estimates of $X^a$ in $C^{1,1}$ open sets in a forthcoming paper
\cite{CKSZ}.

For $a>0$, $X^a$ and $X:=X^1$ are in fact related by a scaling. More
precisely, for $a\in (0, M]$,  $X^a$ has the same distribution as
$\lambda X_{\lambda^{-2}t}$, where $\lambda = a^{\alpha/(\alpha-2)}
\geq M^{\alpha/(\alpha-2)}$. Consequently, if  $u$ is harmonic in an
open set $U$ with respect to  $X^a$, then $v(x):=u(\lambda x)$ is
harmonic in $\lambda^{-1}U$ with respect to $X$. Hence the uniform
Harnack inequality of Proposition \ref{uhp} follows from the Harnack
inequality for $X$. The latter is known, see Theorem 6.7 of
\cite{CK08} or Theorem 4.5 of \cite{SV07}. However the uniform BHP
of Theorem \ref{t:main} can not be obtained by such a scaling
argument from the BHP of $X$. This is because for a $C^{1,1}$ open
set $D$ with the characteristics $(R, \lambda)$, $\lambda^{-1}D$ is,
in general, a $C^{1,1}$ open set with $C^{1,1}$ characteristics
$(R/\lambda, \, \lambda \Lambda)$, which tends to $(0, \infty)$ as
$\lambda \to \infty$.

For each fixed $\alpha_0\in (0, 2)$, when $\alpha$ changes from
$\alpha_0$ to 2, the operator $\Delta +a^\alpha \Delta^{\alpha/2}$
evolves continuously from $\Delta + a^{\alpha_0}
\Delta^{\alpha_0/2}$ to $(1+a^2)\Delta$. So in view of Theorem
\ref{t:main}, it is reasonable to expect that one can get the BHP
for $\Delta +a^\alpha \Delta^{\alpha/2}$  uniformly   both in $a\in
(0, M]$    and  in $\alpha \in [\alpha_0, 2)$. We believe this is
the case and that it can be achieved by carefully keeping track of
all the comparison constants in the arguments of this paper. However
in order to keep our exposition as transparent as possible, we are
content with establishing the results stated in Theorem \ref{t:main}
and leave the details of the proof for the last claim to interested
readers.

Our method of establishing the above BHP is different from those in
\cite{Bo, SW99} for symmetric stable processes and in \cite{KSV} for
more general subordinate Brownian motions. The reason that the
approaches in \cite{Bo, SW99, KSV} do not work well in our setting
lies exactly with the fact that $X^a$ leaves open set $D$ by jumping
out across the boundary $\partial D$ as well as by continuously
exiting $D$ through the boundary of $D$. To circumvent this
difficulty, in this paper we adopt the ideas from \cite{BBC} for the
BHP of censored stable processes, which are further refined in
\cite{G}. That is, we use suitably chosen subharmonic and
superharmonic functions of the process $X^a$ (or equivalently, of
$\Delta + a^{\alpha} \Delta^{\alpha/2}$) to derive some exit
distribution estimates that are needed to establish the BHP.
However, had we done it in this way directly, we would only get the
BHP for $\Delta +a^{\alpha} \Delta^{\alpha/2}$ with $\alpha \in (1,
2)$. The reason is that, when $D=\bH^d_+:=\{x=(x_1,\cdots, x_d)\in
\R^d: \, x_1>0\}$, we need to consider testing functions $w_p(x)=
(x_1\vee 0)^p$ for $p>1$. But for $w_p$ to be
$\Delta^{\alpha/2}$-differentiable in $\bH^d_+$, see \eqref{e:1.1},
one requires $p<\alpha$, which would be impossible when  $\alpha \in
(0, 1]$. To overcome this difficulty, for each $\lambda >0$, we
consider the finite range (or truncated) symmetric $\alpha$-stable
process $\wh Y^\lambda$ obtained from $Y$ by suppressing all its
jumps of size larger than $\lambda$. The infinitesimal generator of
$\wh Y^\lambda$ is
\begin{equation}\label{e:1.2}
\wh \Delta_{d, \lambda}^{\alpha/2} u(x) := \lim_{\eps \downarrow
0}\int_{\{y\in \bR^d: \, \eps <|y-x|<\lambda \}} (u(y)-u(x))
\frac{\sA(d, \alpha)}{|x-y|^{d+\alpha}}\, dy .
\end{equation}
When $\lambda =1$, we will simply denote $\wh \Delta_{d,
1}^{\alpha/2}$ by $\wh \Delta_d^{\alpha/2}$. Then $w_p$ is $\wh
\Delta_d^{\alpha/2}$-differentiable in $\bH^d_+$ for every $p>0$.
Observe that $\wh X^a:=X^0+ a \wh Y^{1/a}$ is a L\'evy process
obtained from $X^a=X^0+aY$ by suppressing all its jumps of size
larger than 1 and that the infinitesimal generator of $\wh X^a$ is
$\Delta + a^\alpha \wh \Delta_d^{\alpha /2}$. From this, we can
obtain suitable exit distribution estimates for the L\'evy process
$\wh X^a$. The desired estimates for $X^a$ can then be obtained from
that for $\wh X^a$ by adding back those jumps of $X^a$ of size
larger than 1. Such an idea has already been used in \cite{CR} to
study Schramm-L\"owner evolutions driven by one-dimensional
symmetric stable processes. We remark that the BHP in Theorem
\ref{t:main} for the case of $a=1$ has also been mentioned in Remark
5.2 of Guan \cite{G}. However, no precise statement (such as the
range of $\alpha$) nor a proof is given in that paper.

The rest of the paper is organized as follows. In Section 2, we
derive estimates on $\wh \Delta_d^{\alpha/2} w_p$. These estimates
are then used in Section 3 to obtain exit distribution (or harmonic
measure) estimates for the finite range process $\wh X^a$ and then
for the desired process $X^a$. In Section 4, we first give the proof
of Proposition \ref{uhp}, and then establish a Carleson estimate for
non-negative harmonic functions of $\Delta+ a^\alpha
\Delta^{\alpha/2}$ in Lipschitz open sets. Then using these results,
the proof of Theorem \ref{t:main} is presented.

Throughout this paper, we use the capital letters $C_1,C_2, \cdots $
to denote constants in the statement of the results, and their
labeling will be fixed. The lowercase constants $c_1, c_2, \cdots$
will denote generic constants used in the proofs, whose exact values
are not important and can change from one appearance to another. The
labeling of the constants $c_1, c_2, \cdots$ starts anew in every
proof. The dependence of the constant $c$ on the dimension $d \ge 1$
and $\alpha \in (0, 2)$ may not be mentioned explicitly. The
constant $M>0$ will be fixed throughout
 this paper. We will use ``$:=$" to denote a definition, which is
read as ``is defined to be". For $a, b\in \bR$, $a\wedge b:=\min
\{a, b\}$ and $a\vee b:=\max\{a, b\}$. For every function $f$, let
$f^+:=f \vee 0$.   We will use $\partial$ to denote a cemetery point
and for every function $f$, we extend its definition to $\partial$
by setting $f(\partial )=0$. We will use $dx$ or $m_{d}(dx)$ to
denote the Lebesgue measure in $\bR^d$. For a Borel set $A\subset
\bR^d$, we also use $|A|$ to denote its Lebesgue measure and ${\rm
diam}(A)$ to denote the diameter of the set $A$.

\section{Truncated fractional Laplacian estimates for power functions}

In this section, we give some estimates which will be used later.
Recall that the fractional Laplacian $\Delta^{\alpha/2}$ and the
truncated fractional Laplacian $\wh \Delta_d^{\alpha/2}:=\wh
\Delta_{d, 1}^{\alpha/2}$ are defined in \eqref{e:1.1} and
\eqref{e:1.2}, respectively.

\begin{lemma}\label{L:2.8}
For $x\in \R^d$ and $p>0$, set $w_p(x): =
 (x_1^+)^{p}$. Then there are constants $R_* \in (0, 1)$, $C_1>
C_2>0$ depending only on $p$, $d$ and $\alpha$ such that for every
$x\in \R^d$ with $x_1 \in (0,R_*]$
\begin{equation}\label{eqn:2.a}
| \wh \Delta^{\alpha/2}_d w_p(x) |  \leq C_1
\qquad \hbox{for } p > \alpha,
\end{equation}
\begin{equation}\label{eqn:2.b_1}
 | \wh \Delta^{\alpha/2}_d w_p(x) |  \leq C_1 \, | \log x_1 |
\qquad \hbox{for } p = \alpha,
\end{equation}
\begin{equation}\label{eqn:2.b}
C_2 x_1^{p-\alpha}  \leq \wh \Delta^{\alpha/2}_d w_p(x)  \leq C_1
x_1^{p-\alpha}  \qquad \hbox{for } \alpha/2 < p < \alpha,
\end{equation}
\begin{equation}\label{eqn:2.d}
- C_1   \leq \wh \Delta^{\alpha/2}_d w_p(x)  \leq - C_2
   \qquad \hbox{for } p=\alpha/2 ,
\end{equation}
and
\begin{equation}\label{eqn:2.c}
-C_1 x_1^{p-\alpha}  \leq \wh \Delta^{\alpha/2}_d w_p(x)  \leq -C_2
x_1^{p-\alpha}  \qquad \hbox{for } 0 < p <\alpha/2.
\end{equation}
\end{lemma}

\pf
First note that using integration by parts and a change of variable,
we get that for $p,x>0$ and $\eps \in (0, 1/(x+1))$,
\begin{eqnarray}
 \int_0^{1-\eps} \frac{z^p-1}{(1-z)^{\alpha+1}}dz
 & =& \frac1{\alpha}
\int_0^{1-\eps}(z^p-1) d(1-z)^{-\alpha} \label{e:FH1}\\
& =& \frac1{\alpha}(z^p-1)(1-z)^{-\alpha} \Big|^{1-\eps}_0
 -\frac{p}{\alpha}
\int_0^{1-\eps} z^{p-1}(1-z)^{-\alpha}dz\nonumber\\
& = & \frac{(1-\eps)^p-1}{\alpha\eps^\alpha}+\frac{1}{\alpha }
 -\frac{p}{\alpha}
\int_{0}^{1-\eps} \frac{z^{p-1}}{(1-z)^{\alpha}}dz \nonumber
\end{eqnarray}
and
\begin{eqnarray}
  && \int_{1+\eps}^{\frac{x+1}{x}} \frac{z^p-1}{(z-1)^{\alpha+1}} dz
  =
- \frac1{\alpha} \int_{1+\eps}^{\frac{x+1}x}
(z^p-1) d(z-1)^{-\alpha} \label{e:FH2}\\
&& =-\frac1{\alpha}(z^p-1)(z-1)^{-\alpha}
 \Big|_{1+\eps}^{\frac{x+1}{x}} +\frac{p}{\alpha}
\int_{1+\eps}^{\frac{x+1}{x}}  z^{p-1}(z-1)^{-\alpha}dz\nonumber\\
&&=\frac{(1+\eps)^p-1}{\alpha\eps^\alpha}+ \frac1{\alpha} x^{\alpha}
-\frac1{\alpha} (x+1)^px^{\alpha-p}+\frac{p}{\alpha}
\int_{\frac{x}{x+1}}^{\frac{1}{\eps+1}}
z^{\alpha-p-1}(1-z)^{-\alpha}dz . \nonumber
\end{eqnarray}
For $p>0$ and $x \in (0,1)$,   by a change of variable
\begin{align*}
& \wh \Delta^{\alpha/2}_1 w_p(x)
 = \sA(1, \alpha)\,\lim_{\eps
 \downarrow 0}\int_\R \frac{ w_p(y)-w_p(x)}{|x-y|^{1+\alpha}} \,
\1_{\{\eps<|y-x|\leq 1\}}\, dy \\
&= \sA(1, \alpha)\,\lim_{\eps \downarrow 0}\int_{0}^{x+1} \frac{
 y^p-x^p}{|x-y|^{1+\alpha}} \, \1_{\{ |y-x|>\eps\}}\, dy -
\sA(1, \alpha)x^p\int_{x-1}^{0} \frac{dy}{|x-y|^{1+\alpha}}\\
&= \sA(1, \alpha)\,x^{p-\alpha}\, \lim_{\eps \downarrow 0}
\int_{0}^{\frac{x+1}{x}} \frac{ z^p-1}{|z-1|^{1+\alpha}}  \,
\1_{\{|z-1|>\eps/x\}}\, dz -
\sA(1, \alpha)x^{p} \int_{x-1}^0 (x-y)^{-1-\alpha} dy\\
&= \sA(1, \alpha)\,x^{p-\alpha} \, \lim_{\eps \downarrow
0}\left(\int_{0}^{1-\eps} \frac{
z^p-1}{(1-z)^{1+\alpha}}dz+\int_{1+\eps}^{\frac{x+1}{x}} \frac{
z^p-1}{(z-1)^{1+\alpha}}dz\right)- \frac{\sA(1, \alpha)}{\alpha}
(x^{p-\alpha}-x^p).
\end{align*}
So we have by   \eqref{e:FH1}-\eqref{e:FH2} that for $p>0$ and $x
\in (0, 1)$,
\begin{eqnarray}
&&\wh \Delta^{\alpha/2}_1 w_p(x) \label{e:C1}\\
&&= \frac{\sA(1, \alpha)}{\alpha} \,x^{p-\alpha} \lim_{\eps
\downarrow 0} \left(1+ \frac{(1-\eps)^p+ (1+\eps)^p-2}{\eps^\alpha}
\right) - \frac{\sA(1, \alpha)}{\alpha}
(x^{p-\alpha}-x^p)+\frac{\sA(1, \alpha)}{\alpha} x^{p} \nonumber\\
&& \,\,\,\,\,\,-\frac{\sA(1, \alpha)}{\alpha} (x+1)^p +\frac{
\sA(1, \alpha)\, p}{\alpha} x^{p-\alpha} \left(
\int_{\frac{x}{x+1}}^{1} \frac{
z^{\alpha-p-1}-z^{p-1}}{(1-z)^{\alpha}}dz-\int_{0}^{\frac{x}{x+1}}
\frac{z^{p-1}}{(1-z)^{\alpha}}dy\right) \nonumber\\
&&=\frac{\sA(1, \alpha)}{\alpha} \left(2 x^p
-(x+1)^p+{p}x^{p-\alpha} \Big( \int_{\frac{x}{x+1}}^{1} \frac{
z^{\alpha-p-1}-z^{p-1}}{(1-z)^{\alpha}}dz-\int_{0}^{\frac{x}{x+1}}
\frac{ z^{p-1}}{(1-z)^{\alpha}}dz\Big)\right). \nonumber
\end{eqnarray}
 Note that for $p>\alpha$,
$$
\sup_{x\in (0, 1]} x^{p-\alpha} \left| \int_{\frac{x}{x+1}}^{1}
\frac{z^{\alpha-p-1}-z^{p-1}}{(1-z)^{\alpha}}dz-\int_{0}^{\frac{x}{x+1}}
\frac{ z^{p-1}}{(1-z)^{\alpha}}dz \right| <\infty.
$$
So when $p> \alpha$,
\begin{equation}\label{e:C3}
\sup_{x\in (0, 1)} |\wh \Delta^{\alpha/2}_1 w_p(x)|<\infty  .
\end{equation}
When $p= \alpha$, there exists an $r_*>0$ such that  for $0<x<r_*$
\begin{equation}\label{e:C2_1}
\left| \int_{\frac{x}{x+1}}^{1} \frac{
z^{-1}-z^{{\alpha}-1}}{(1-z)^{\alpha}}dz-\int_{0}^{\frac{x}{x+1}}
\frac{ z^{{\alpha}-1}}{(1-z)^{\alpha}}dz \right| \le \int_{\frac{x}{x+1}}^{1} \frac{
z^{-1}}{(1-z)^{\alpha}}dz \le (1+r_*)^\alpha \log ((1+r_*)/x).
\end{equation}
It is easy to see that
\begin{equation}\label{e:C2_2}
\sup_{x\in [r_*, 1]}\left| \int_{\frac{x}{x+1}}^{1} \frac{
z^{-1}-z^{{\alpha}-1}}{(1-z)^{\alpha}}dz-\int_{0}^{\frac{x}{x+1}}
\frac{ z^{{\alpha}-1}}{(1-z)^{\alpha}}dz \right| <\infty.
\end{equation}
On the other hand, when $p\in (0, \alpha)$,
\begin{equation}\label{e:C2}
\sup_{x\in (0, 1]}\left| \int_{\frac{x}{x+1}}^{1} \frac{
z^{\alpha-p-1}-z^{p-1}}{(1-z)^{\alpha}}dz-\int_{0}^{\frac{x}{x+1}}
\frac{ z^{p-1}}{(1-z)^{\alpha}}dz \right| <\infty.
\end{equation}
As
$$
\lim_{x\to 0+} \left(\int_{\frac{x}{x+1}}^{1} \frac{
z^{\alpha-p-1}-z^{p-1}}{(1-z)^{\alpha}}dz-\int_{0}^{\frac{x}{x+1}}
\frac{ z^{p-1}}{(1-z)^{\alpha}}dz \right) \
\begin{cases} \ >0 \qquad &\hbox{if } p \in (\alpha/2, \alpha) \\
 \ <0 &\hbox{if } p\in (0, \alpha/2)
\end{cases}
$$
while for $p=\alpha/2$,
$$ \int_{\frac{x}{x+1}}^{1} \frac{
z^{\alpha-p-1}-z^{p-1}}{(1-z)^{\alpha}}dz-\int_{0}^{\frac{x}{x+1}}
\frac{ z^{p-1}}{(1-z)^{\alpha}}dz =-\int_0^{\frac{x}{x+1}} \frac{z^{\alpha/2-1} }{(1-z)^\alpha} dz,
$$
we conclude from \eqref{e:C1}-\eqref{e:C2} that   there are
constants $r_1\in (0, 1)$ and $C_1> C_2>0$ depending on $p$ and
$\alpha$ so that
  when $p= \alpha$,
\begin{equation}\label{e:C4_1}
  |\wh \Delta^{\alpha/2}_1 w_p(x) < C_1\,
|\log x| \ \hbox{ for } x\in (0, r_1]  \quad \hbox{and} \quad
\sup_{x\in (r_1, 1)} |\wh \Delta^{\alpha/2}_1 w_p(x)|<\infty,
\end{equation}
when $p=\alpha/2$,
\begin{equation}
\label{e:C4_2}
 -C_1 \leq \wh \Delta^{\alpha/2}_1 w_{p} (x) < -C_2
  \ \hbox{ for } x\in (0, r_1]  \quad \hbox{and} \quad
\sup_{x\in (r_1, 1)} |\wh \Delta^{\alpha/2}_1 w_p(x)|<\infty,
\end{equation}
when $p\in (\alpha/2, \alpha)$,
\begin{equation}\label{e:C4}
C_2 \, x^{p-\alpha}< \wh \Delta^{\alpha/2}_1 w_p(x) < C_1\,
x^{p-\alpha} \ \hbox{ for } x\in (0, r_1]  \quad \hbox{and} \quad
\sup_{x\in (r_1, 1)} |\wh \Delta^{\alpha/2}_1 w_p(x)|<\infty,
\end{equation}
and for $p\in (0, \alpha/2)$,
\begin{equation}\label{e:C5}
-C_1 \, x^{p-\alpha}< \wh \Delta^{\alpha/2}_1 w_p(x) < -C_2\,
x^{p-\alpha} \  \hbox{ for } x\in (0, r_1] \quad \hbox{and} \quad
\sup_{x\in (r_1, 1)} |\wh \Delta^{\alpha/2}_1 w_p(x)|<\infty.
\end{equation}

On the other hand, for
$x \ge 1$,
\begin{align*}
\wh \Delta^{\alpha/2}_1 w_p(x) =& \sA(1, \alpha)\,\lim_{\eps
\downarrow 0}\int_{x-1}^{x+1} \frac{
w_p(y)-w_p(x)}{|x-y|^{1+\alpha}}\, \1_{\{x-y|>\eps\}} \, dy\\
=& \sA(1, \alpha)\,x^{p-\alpha}\lim_{\eps \downarrow
0}\int_{\frac{x-1}{x}}^{\frac{x+1}{x}} \frac{
y^p-1}{|y-1|^{1+\alpha}} \, \1_{\{|y-1|>\eps\}} \, dy\\
=& \sA(1, \alpha)\,x^{p-\alpha}\lim_{\eps \downarrow
0}\left(\int_{\frac{x-1}{x}}^{1-\eps} \frac{
y^p-1}{(1-y)^{1+\alpha}}dy+\int_{1+\eps}^{\frac{x+1}{x}} \frac{
y^p-1}{(y-1)^{1+\alpha}}dy\right) \\
=& \sA(1, \alpha)\,x^{p-\alpha} \int_0^{1/x}  \frac{
(1+u)^p+(1-u)^p-2}{u^{1+\alpha}}du .
\end{align*}
Note
the above integrand
$$
\frac{(1+u)^p+(1-u)^p-2}{u^{1+\alpha}}
$$
is of the order $u^{1-\alpha}$ near zero.
 So for $p>0$ and $\alpha \in (0, 2)$, there is a constant
$c_0=c_0(p, \alpha)>0$ so that
\begin{equation}\label{e:C6}
|\wh \Delta^{\alpha/2}_1 w_p(x)| \leq c_0 x^{p-2} \qquad
\hbox{for } x\ge1.
\end{equation}
With $r_1\in (0, 1)$ as in
\eqref{e:C4_1}-\eqref{e:C5}, the above
inequality in fact holds for $x>r_1$.

The estimates \eqref{e:C3}-\eqref{e:C5} prove the Lemma in dimension
$d=1$. Now we consider the case $d\ge 2$. For each fixed $x \in
\R^d$, we use the spherical coordinates
$$
(y_1,\dots,y_d) :=x+(r\cos \theta_1, r\sin \theta_1  \cos \theta_2,
\dots,r \sin \theta_1\dots\cos \theta_{d-1}, r \sin \theta_1\dots
\sin {\theta_{d-1}})
$$
where $r\ge 0$, $0 \le \theta_1, \dots, \theta_{d-2} < \pi$ and
$0\le \theta_{d-1} < 2 \pi$. Let
$$
\phi(\wh \theta):=\phi(\theta_1, \dots, \theta_{d-2}):=\sin^{d-2}
\theta_1\sin^{d-3} \theta_2\dots \sin \theta_{d-2}.
$$
Then for $x\in \R^d$ with $x_1 >0$ we have
\begin{align*}
&\lim_{\eps \downarrow 0}\int_{\{y\in
\bR^d: \, 1>|y-x|>\eps\}} ((y_1^+)^p-x_1^p)\frac{dy}{|x-y|^{d+\alpha}}\\
= &\lim_{\varepsilon\downarrow 0} \int_0^{\pi} d\theta_1\cdots
\int_0^\pi d\theta_{d-2}\int_0^{2\pi} \phi(\wh \theta) d\theta_{d-1}
\int_{\eps}^1 \frac{((r\cos\theta_1 +
x_1)^+)^p- x_1^p}{r^{d+\alpha}} r^{d-1}dr\\
=&\lim_{\varepsilon\downarrow 0} \int_0^{\pi/2}
d\theta_1\int_0^{\pi} d\theta_2\cdots \int_0^\pi
d\theta_{d-2}\int_0^{2\pi} \phi(\wh \theta) d\theta_{d-1}
(\cos\theta_1)^p\int_{\eps}^1 \frac{((r + \frac{x_1}{\cos
\theta_1})^+)^p- (\frac{x_1}{\cos
\theta_1})^p}{r^{1+\alpha}}dr\\
&+\lim_{\varepsilon\downarrow 0} \int_{\pi/2}^{\pi}
d\theta_1\int_0^{\pi} d\theta_2\cdots \int_0^\pi
d\theta_{d-2}\int_0^{2\pi} \phi(\wh \theta) d\theta_{d-1}
(-\cos\theta_1)^p\int_{\eps}^1 \frac{((-r - \frac{x_1}{\cos
\theta_1})^+)^p- (-\frac{x_1}{\cos \theta_1})^p}{r^{1+\alpha}}dr.
\end{align*}
By the change of variable $r=t-x_1/\cos \theta_1$ for $\theta \in
[0, \pi/2)$ and $r=-t-x_1/\cos \theta_1=-t+x_1/\cos (\pi- \theta_1)
$ for $\theta \in (\pi/2, \pi]$, we get
\begin{align*}
&\lim_{\eps \downarrow 0}\int_{\{y\in
\bR^d: \, 1>|y-x|>\eps\}} ((y_1^+)^p-x_1^p)  \frac{dy}{|x-y|^{d+\alpha}}\\
=&\lim_{\varepsilon\downarrow 0} \int_0^{\pi/2}
d\theta_1\int_0^{\pi} d\theta_2\cdots \int_0^\pi
d\theta_{d-2}\int_0^{2\pi} \phi(\wh \theta) d\theta_{d-1}
(\cos\theta_1)^p\int_{\eps+\frac{x_1}
{\cos\theta_1}}^{1+\frac{x_1}{\cos \theta_1}} \frac{(t^+)^p-
(\frac{x_1}
{\cos \theta_1})^p}{|t-x_1/\cos \theta_1|^{1+\alpha}}dt\\
&+\lim_{\varepsilon\downarrow 0} \int_{\pi/2}^{\pi}
d\theta_1\int_0^{\pi} d\theta_2\cdots \int_0^\pi
d\theta_{d-2}\int_0^{2\pi} \phi(\wh \theta) d\theta_{d-1}
(\cos(\pi-\theta_1))^p \int^{-\eps+\frac{x_1}{\cos(\pi-\theta_1)}}_{
-1+\frac{x_1}{\cos (\pi-\theta_1)}} \frac{(t^+)^p- (\frac{x_1}{\cos
(\pi-\theta_1)})^p}{|t-x_1/\cos(\pi- \theta_1)|^{1+\alpha}}dt\\
=&\lim_{\varepsilon\downarrow 0} \int_0^{\pi/2}
d\theta_1\int_0^{\pi} d\theta_2\cdots \int_0^\pi
d\theta_{d-2}\int_0^{2\pi} \phi(\wh \theta) d\theta_{d-1}
(\cos\theta_1)^p\int_{\eps+\frac{x_1}
{\cos\theta_1}}^{1+\frac{x_1}{\cos \theta_1}} \frac{(t^+)^p-
(\frac{x_1}
{\cos \theta_1})^p}{|t-x_1/\cos \theta_1|^{1+\alpha}}dt\\
&+\lim_{\varepsilon\downarrow 0} \int_{0}^{\pi/2}
d\theta_1\int_0^{\pi} d\theta_2\cdots \int_0^\pi
d\theta_{d-2}\int_0^{2\pi} \phi(\wh \theta) d\theta_{d-1}
(\cos\theta_1)^p\int^{-\eps+\frac{x_1} {\cos\theta_1}}_{-1+\frac{x_1}
{\cos \theta_1}} \frac{(t^+)^p- (\frac{x_1}{\cos \theta_1})^p}
{|t-x_1/\cos\theta_1|^{1+\alpha}}dt\\
=&\int_{0}^{\pi/2} d\theta_1\int_0^{\pi} d\theta_2\cdots \int_0^\pi
d\theta_{d-2}\int_0^{2\pi} \phi(\wh \theta) d\theta_{d-1}
(\cos\theta_1)^p\left(\lim_{\varepsilon\downarrow 0} \int_{\{t\in
\bR: \, 1>|t-\frac{x_1}{\cos \theta_1}|>\eps\}}  \frac{(t^+)^p-
(\frac{x_1}{\cos
\theta_1})^p}{|t-x_1/\cos\theta_1|^{1+\alpha}}dt\right).
\end{align*}
Therefore we have
\begin{eqnarray*}
&&  \wh \Delta^{\alpha/2}_d w_p(x) \\
&=& \frac{\sA(d,\alpha)}{\sA(1,\alpha)} \int_{0}^{\pi/2}
d\theta_1\int_0^{\pi} d\theta_2\cdots \int_0^\pi
d\theta_{d-2}\int_0^{2\pi}
\phi(\wh \theta) d\theta_{d-1} (\cos\theta_1)^p \wh
\Delta^{\alpha/2}_1 w_p \left(\frac{x_1}{\cos \theta_1}\right) \\
&=& \frac{\sA(d,\alpha)}{\sA(1,\alpha)} \int_{0}^{\arccos (x_1/r_1)}
d\theta_1\int_0^{\pi} d\theta_2\cdots \int_0^\pi
d\theta_{d-2}\int_0^{2\pi}
\phi(\wh \theta) d\theta_{d-1} (\cos\theta_1)^p \wh \Delta^{\alpha/2}_1
w_p \left(\frac{x_1}{\cos \theta_1}\right) \\
&& +\frac{\sA(d,\alpha)}{\sA(1,\alpha)} \int_{\arccos
(x_1/r_1)}^{\pi/2} d\theta_1\int_0^{\pi} d\theta_2\cdots \int_0^\pi
d\theta_{d-2}\int_0^{2\pi} \phi(\wh \theta) d\theta_{d-1}
(\cos\theta_1)^p \wh \Delta^{\alpha/2}_1 w_p \left(\frac{x_1}{\cos
\theta_1}\right).
\end{eqnarray*}
The conclusion \eqref{eqn:2.a}-\eqref{eqn:2.c} now follow
immediately from the above equality and the estimates
\eqref{e:C3}-\eqref{e:C6}, where we use \eqref{e:C6} to bound the
second integral above by $c_1 x_1^3/r_1^3$ for some positive
constant $c_1$. \qed

\begin{remark}\label{R:2.2}\rm
A careful evaluation of \eqref{e:C1} in fact shows that $\lim_{x\to
0^+}\wh \Delta_1^{\alpha/2}w_p(x)=c\not= 0$ when $p>\alpha$. At the
first glance, this   may look surprising, as in the Brownian motion
case (which corresponds to $\alpha =2$), $\Delta x_1^p =
p(p-1)x^{p-2}_1$. The bound in \eqref{eqn:2.a} is due to the
non-local nature of the operator $\wh \Delta_d^{\alpha/2}$ for
$\alpha \in (0, 2)$. However a more careful analysis of \eqref{e:C1}
reveals that for $p>0$,
$$
\wh \Delta_1^{\alpha/2} w_p(x) \asymp
(2-\alpha)\left((p-\alpha)^{-1}-1\right)+ p(p-1)x^{p-\alpha} \qquad
\hbox{for }   x\in (0, r_1)
$$
as $\alpha \uparrow 2$. It is not difficult to see that as $\alpha
\uparrow 2$, $\wh \Delta_1^{\alpha/2} w_p(x) $ converges to
 $\Delta w_p(x) $. \qed
\end{remark}

Recall that for $\lambda >0$, the operator $\wh
\Delta_{d,\lambda}^{\alpha/2}$ is defined by \eqref{e:1.2}.
 Note that
\begin{align}
\wh \Delta_{d,\lambda}^{\alpha/2} u(x) =\lambda^{-\alpha}(\widehat
\Delta_{d}^{\alpha/2}u(\lambda\cdot)) (\lambda^{-1} x). \label{e:L1}
\end{align}
Thus, from Lemma \ref{L:2.8} and \eqref{e:L1}, we get the following
corollary.

\begin{cor}\label{C:2.8}
For $x\in \R^d$ and $p>0$, set $w_p(x): =
 (x_1^+)^{p}$. Then there are constants $R_* \in (0, 1/2)$, $C_1>
C_2>0$ depending only on $p$, $d$ and $\alpha$ such that for every
$\lambda >0$ and $x\in \R^d$ with
$x_1 \in (0,\lambda R_* )$,
\begin{equation}\label{eqn:2.a1}
|\wh \Delta^{\alpha/2}_{d, \lambda} w_p(x)|  \leq C_1
\lambda^{p-\alpha} \qquad \hbox{for } p > \alpha ,
\end{equation}
\begin{equation}\label{eqn:2.d1}
 |\wh \Delta^{\alpha/2}_{d, \lambda} w_p(x)|
\le C_1  |\log  (x_1/\lambda)|, \qquad \hbox{for } p =
\alpha,
\end{equation}
\begin{equation}\label{eqn:2.b1}
C_2 x_1^{p-\alpha}   \leq \wh \Delta^{\alpha/2}_{d, \lambda}
w_p(x) \leq C_1 x_1^{p-\alpha}  \qquad \hbox{for } \alpha/2 < p <
\alpha,
\end{equation}
\begin{equation}
\label{eqn:2.e1}
-C_1  \lambda^{-\alpha/2} \leq \wh \Delta^{\alpha/2}_{d, \lambda}
w_p(x) \leq -C_2 \lambda^{-\alpha/2}  \qquad \hbox{for } p=\alpha/2 ,
\end{equation}
and
\begin{equation}\label{eqn:2.c1}
-C_1 x_1^{p-\alpha}   \leq \wh \Delta^{\alpha/2}_{d, \lambda}
w_p(x)  \leq -C_2 x_1^{p-\alpha}  \qquad \hbox{for } 0 < p <
\alpha/2.
\end{equation}
\end{cor}

\section{Estimates on harmonic measures}

Recall that for any open set $U\subset \bR^d$,
$\tau^a_U=\inf\{t>0: \, X^a_t\notin U\}$ is  the first exit time
from $U$ by $X^a$.

\begin{lemma}\label{L:2.00}
For every $b \in (0,\infty)$, there
 exist $C_3=C_3(M, b)>0$ and $C_4=C_4(M, b)>0$ such that for every
$x_0 \in \bR^d$, $a \in (0,M]$ and $r \in (0, b]$,
 \bee \label{e:ext}
C_3 r^2 \, \le \, \E_{x_0}\left[\tau^a_{B(x_0,r)}\right]\, \le\,
C_4\,r^2.
 \eee
\end{lemma}
\pf See Lemmas 2.2 and 2.3 in \cite{SV05} or Lemmas 2.3 and 2.4 in
\cite{CK08} for a proof. \qed

In the remainder of this section, we assume $D$ is a $C^{1,1}$ open
set with characteristics $(R, \Lambda)$. Recall that we are always
assuming that $R \le 1$ and $\Lambda \ge 1$. For notational
convenience, throughout the rest of this section, we put
$$
r_0= \frac{R}{4\sqrt{1+\Lambda^2}}.
$$
Define
$$
\rho_Q (x) := x_d -  \phi_Q (\wt x),
$$
where $(\wt x, x_d)$ is the coordinates of $x$ in $CS_Q$. Note that
for every $Q \in \partial D$ and $ x \in B(Q, R)\cap D$ we have
\begin{equation}\label{e:d_com}
(1+\Lambda^2)^{-1/2} \,\rho_Q (x) \,\le\, \delta_D(x)  \,\le\,
\rho_Q(x).
\end{equation}

Recall that  $R_*$ is the constant in Lemma \ref{L:2.8}.

\begin{lemma}\label{L:Main}
Fix $Q \in \partial D$ and the coordinate system $CS_Q$ so that
$$
B(Q, R)\cap D = \left\{ y=(\wt y, \, y_d) \in B(0, R) \mbox{ in }
CS_Q: y_d > \phi (\wt y) \right\}.
$$
For $p>\alpha/2$, let
$$
h_p(y):=\left(\rho_Q (y)\right)^p{\bf 1}_{D\cap B(Q, 4r_0)}(y).
$$
Then there exist $C_i=C_i(\alpha, p, \Lambda, R)>0$, $i=5,6,7$,
independent of the choice of the point $Q \in \partial D$ such that

\begin{description}
\item {\rm (i)} in the case $\frac{\alpha}2< p<\alpha$, for
all $x\in D$ such that
$\rho_Q(x)<r_0 \wedge R_*$ and
$|\widetilde{x}| <r_0$, we have
 \bee\label{e:h2}
C_6 \left(\rho_Q (x)\right)^{p-\alpha}  \leq \wh \Delta_d^{\alpha/2}
h_p(x)  \le C_5 \left(\rho_Q (x)\right)^{p-\alpha};
 \eee
\item{\rm (ii)} in the case $ p> {\alpha} $, for
all $x\in D$ such that
$\rho_Q(x)<r_0 \wedge R_*$ and
$|\widetilde{x}| <r_0$, we have
\bee
\label{e:h0} |\wh \Delta_d^{\alpha/2} h_p(x)|  \le C_7 ;
\eee
\item{\rm (iii)} in the case $ p= {\alpha} $, for
all $x\in D$ such that
$\rho_Q(x)<r_0 \wedge R_*$ and
$|\widetilde{x}| <r_0$, we have
\bee
\label{e:h0_1} |\wh \Delta_d^{\alpha/2} h_p(x)|  \le C_7 |\log \left(\rho_Q (x)\right)|\, .
\eee
\end{description}
\end{lemma}

\pf In this proof our coordinate system is always $CS_Q$. Fix
$x=(\wt{x}, x_d)\in D$ such that $\rho_Q(x)<r_0 \wedge R_*$ and
$|\wt{x}|<r_0$, and choose a point $x_0\in\partial D$ satisfying
$\widetilde{x}=\widetilde{x_0}$. Denote by $\overrightarrow{n}(x_0)$
the inward   unit normal  vector at $x_0$ for $\partial D$ and set
$\Phi(y)=\langle y-x_0,\overrightarrow{n}(x_0)\rangle$ for $y\in
\bR^d$. Then $\Pi:=\{y:\Phi(y)=0\}$ is the hyperplane tangent to
$\partial D$ at the point $x_0$. The function
$\Gamma^*:\bR^{d-1}\rightarrow \bR$ describing the plane $\Pi$ is
given by $\Gamma^*(\tilde{y})=\phi_Q(\tilde{x}_0)
+\nabla\phi_Q(\tilde{x}_0)(\tilde{y}-\tilde{x}_0)$, and it holds
that
 $ \big\langle
 \left(\widetilde{y},\Gamma^*(\widetilde{y})\right)-x_0, \,
\overrightarrow{n}(x_0)\big\rangle=0. $ We also let
\begin{eqnarray*}
A&:=&\big\{y: \Gamma^*(\widetilde{y})<y_d<\phi_Q(\widetilde{y})
\hbox{ and } |\widetilde{y}-\widetilde{x}|<r_0 \big\} \bigcup
\big\{y: \Gamma^*(\widetilde{y})>y_d>\phi_Q(\widetilde{y}) \hbox{
and }
|\widetilde{y}-\widetilde{x}|<r_0 \big\}, \\
E &:=& \big\{y\in D \setminus A:   \
|\widetilde{y}-\widetilde{x}|<r_0 \hbox{ and }  \rho_Q (y) <
r_0(2+\Lambda) \big\}.
\end{eqnarray*}
Note that, if $|x-y| < r_0$ and $y \in D$,
$$
\rho_Q(y)\le |y_d-x_d|+|x_d- \phi_Q(\wt x)|+|\phi_Q(\wt
y)-\phi_Q(\wt x)| < r_0(2+\Lambda).
$$
On the other hand if $|\widetilde{y}- \widetilde{x}|<r_0$ and
$\rho_Q (y) < r_0 (2+\Lambda)$, then
$$
|y|^2=|\widetilde{y}|^2+|y_d|^2 \le
(2r_0)^2+(r_0(2+\Lambda)+|\phi_Q(\wt y)|)^2 \le
\frac{4+(2+3\Lambda)^2}{64(1+\Lambda^2)} R^2 <R^2.
$$
Consequently, we have
 \bee\label{e:BE}  D \cap B(x, r_0) \subset D \cap
\big\{y:\, |\widetilde{y}- \widetilde{x}|<r_0 \hbox{ and } \rho_Q (y) <
r_0 (2+\Lambda) \big\} \subset D \cap B(0,R). \eee

Let
$\overline{h}(y):=\overline{h}_x(y):=(y_d-\Gamma^*(\widetilde{y}))^+$
for $y\in \bR^d$. Since
$\nabla\phi_Q(\widetilde{x})=\nabla\Gamma^*(\widetilde{x})$, by the
mean value theorem and the $C^{1,1}$ condition on $\phi_Q$,
\begin{align}\label{hhhh}
&|\overline{h}({y})-\rho_Q({y})|\, =\, |\phi_Q(\widetilde{y}) -
\Gamma^*(\widetilde{y})| \\
=\,&|\phi_Q(\widetilde{y}) - \phi_Q(\widetilde{x}) -
\nabla\phi_Q(\widetilde{x})\cdot (\widetilde{y}-\widetilde{x})|\,
\leq\, \Lambda|\widetilde{y}-\widetilde{x}|^{2},\qquad  y\in
E.\nonumber
\end{align}
For $y\in \bR^d$, define
$\delta_{_{\Pi}}(y):={ \rm dist}(y,\Pi)$
and
${D}_{\Gamma^*}= \left\{y\in \bR^d:y_d>\Gamma^*(\widetilde{y})\right\}$.
Let
$$
 b_x:=\left(1+|\nabla\phi_Q({\widetilde{x}})|^2\right)^{1/2}
\qquad
\hbox{and} \qquad h_{x,p}(y):=(\overline{h}(y))^p
\quad \hbox{for }
p > \alpha/2.
$$
Note that
$1\le  b_x \le \sqrt{1+\Lambda^2}$
and $h_{x,p}(x)=h_p(x)$.

Recall that $R_*$ and $C_1>C_2>0$ are the constants in Lemma
\ref{L:2.8}. Since $\overline{h}(y)=b_x\delta_{_{\Pi}}(y)$ on
${D}_{\Gamma^*}$, by Lemma \ref{L:2.8}, it holds that for $y\in
{D}_{\Gamma^*}$ and $\delta_{_{\Pi}}(y)< R_*$,
\begin{align}
C_2 \,b_x^p (\delta_{_{\Pi}}(y))^{p-\alpha} \leq \wh
\Delta_d^{\alpha/2} h_{x,p}(y) \,=\, b_x^p \wh \Delta_d^{\alpha/2}
(\delta_{_{\Pi}}(y))^p \,\le\,C_1 \,b_x^p
(\delta_{_{\Pi}}(y))^{p-\alpha}   \quad \text{when }\alpha/2 < p
<\alpha,  \label{e:ee1}
\end{align}
\begin{align}
|\wh \Delta_d^{\alpha/2} h_{x,p}(y)|\,=\,b_x^p | \wh
\Delta_d^{\alpha/2} (\delta_{_{\Pi}}(y))^p |\,\le\,C_1
\,b_x^p\,\le\,C_1(1+\Lambda^2)^{p/2}   \qquad \text{when } p>
\alpha ,  \label{e:ee2}
\end{align}
\begin{eqnarray}\label{e:ee0}
|\wh \Delta_d^{\alpha/2} h_{x,p}(y)|&=&b_x^p | \wh
\Delta_d^{\alpha/2} (\delta_{_{\Pi}}(y))^p |\,\le\,C_1
\,b_x^p|\log(\delta_{_{\Pi}}(y))| \\
&\le&C_1(1+\Lambda^2)^{p/2}|\log(\delta_{_{\Pi}}(y))|
\hskip 0.8truein \text{when } p= \alpha.  \nonumber
\end{eqnarray}
Note that
$b_x\delta_{_{\Pi}}(x)=\rho_Q (x)$.
 Applying \eqref{e:d_com} and
\eqref{e:ee1} to the point $x$ gives that, for $\alpha/2<p<\alpha$
\begin{equation}\label{e:ee1x}
C_2\rho_Q(x)^{p-\alpha} \le C_2 b_x^{\alpha}\rho_Q(x)^{p-\alpha} \le
\widehat{\Delta}^{\alpha/2}_d h_{x,p}(x) \le C_1
b_x^{\alpha}\rho_Q(x)^{p-\alpha}\le
C_1(1+\Lambda^2)^{\alpha/2}\rho_Q(x)^{p-\alpha}.
\end{equation}
Note that by \eqref{e:BE},
\begin{eqnarray}
&& | \wh \Delta_d^{\alpha/2} (h_p-h_{x,p})(x)| \label{e:I}  \\
&=& \sA (d, \alpha)  \bigg|\lim_{\varepsilon\downarrow 0}\int_{
\{1\geq|y-x|>\varepsilon\}}\frac{(h_p(y)-h_{p,x}(y))}{|x-y|^{d+\alpha}}\
dy\bigg|\nonumber\\
&\leq & \sA (d, \alpha)  \bigg| \int_{
\{1\geq |y-x|>r_0\}}\frac{(h_p(y)-h_{p,x}(y))}{|x-y|^{d+\alpha}}\
dy\bigg|\nonumber\\
&&+\sA (d, \alpha)   \lim_{\varepsilon\downarrow 0}\int_{
\{r_0\geq |y-x|>\varepsilon\}}\frac{|h_p(y)-h_{p,x}(y)|}{|x-y|^{d+\alpha}}\
dy \nonumber\\
&\leq & \sA (d, \alpha)  \bigg| \int_{
\{1\geq |y-x|>r_0\}}\frac{(h_p(y)-h_{p,x}(y))}{|x-y|^{d+\alpha}}\
dy\bigg|\nonumber\\
&&+\sA (d, \alpha)\int_{A} \frac{h_p(y)+h_{p,x}(y)}
{|x-y|^{d+\alpha}}\ dy +\sA (d, \alpha) \int_{E}
\frac{|h_p(y)-h_{p,x}(y)|}{|x-y|^{d+\alpha}}  \nonumber\\
&=:&I_1+I_2+I_3.\nonumber
\end{eqnarray}
We claim that, if $p >\alpha/2$, then
 \bee\label{e:I1234}
I_1+I_2+I_3 \leq c_0
 \eee for some constant
$c_0=c_0(\alpha,p,\Lambda, R)$.
Together with \eqref{e:ee2}--\eqref{e:I} this
will establish the desired estimates \eqref{e:h2}--\eqref{e:h0_1}
with constants depending on $\alpha$, $p$, $\Lambda$ and $R$.

Clearly $I_1$ is bounded by some positive constant.

For $y \in A$, we have
\begin{align}
&|h_{x,p}(y)|+|h_{p}(y)|  \ \le \ |y_d-\Gamma^*(\widetilde{y})|^p
+|y_d-\phi_Q(\widetilde{y})|^p \ \le \  2|\phi_Q(\widetilde{y})-
\Gamma^*(\widetilde{y})|^p \label{e:dfe2} \\
\ \le\ 2 & |\phi_Q(\widetilde{y})-\phi_Q(\widetilde{x})-
\nabla\phi_Q(\widetilde{x})\cdot (\widetilde{y}-\widetilde{x})|^p \
\le \  2 \Lambda^{p}|\widetilde{y}-\widetilde{x}|^{2p}
.\nonumber
\end{align}
Furthermore, since, on $\{|\widetilde{y}-\widetilde{x}|=r\le r_0\}$,
$|\phi_Q(\widetilde{y})-\Gamma^*(\widetilde{y})| \le
\Lambda|\widetilde{y}-\widetilde{x}|^{2}=\Lambda r^2$,
\begin{align*}
m_{d-1} \left( \left\{ y:\, | \widetilde{y} -
\widetilde{x}|=r,\Gamma^*(\widetilde{y})<y_d<\phi_Q(\widetilde{y})
\mbox{ or } \Gamma^*(\widetilde{y})>y_d> \phi_Q(\widetilde{y})
\right\} \right)\le c_1 r^{d}
\end{align*}
 for some constant $c_1>0$ if  $r \le r_0$.
This together with \eqref{e:dfe2} yields that
\begin{eqnarray*}
I_2 &\leq& \sA (d, \alpha)  \int_{ 0}^{r_0}\int_{|\wt y-\wt
x|=r}{\bf
 1}_A(y) \frac{|h_{x,p}(y)|+|h_{p}(y)|}{|\wt y-\wt
x|^{d+\alpha}}
\ m_{d-1}(d y) dr\nonumber\\
&\leq & c_2\int_{ 0}^{ r_0} r^{-d+2p-\alpha} m_{d-1}(\{ y \in A :\,
| \widetilde{y} - \widetilde{x}|=r  \}) dr \\
&\le&  c_1c_2\int_{ 0}^{ r_0} r^{2p-\alpha} dr \le  c_3\, .
\end{eqnarray*}

Note that for $y\in E$
\begin{align}
|h_p(y)-h_{x,p}(y)| = & | (\overline{h}(y))^{p}
-(\rho_Q(y))^{p}| \le c_4  (\overline{h}(y))^{(p-1)_-}|
\overline{h}(y) -\rho_Q(y)|, \label{e:KG}
\end{align}
where $(p-1)_-:=(p-1) \wedge 0$.
In the last inequality above, we have used the inequalities
\begin{align*}
|b^{p}-a^{p}|\leq b^{p-1}|b-a| \qquad \hbox{for } a , b>0, \ 0 <p
\le 1
\end{align*}
and
\begin{align*}
|b^{p}-a^{p}|\leq (p+1) |b-a| \qquad \hbox{for } a,b \in (0, 1),  \
p >1 .
\end{align*}
For $y=(\wt y, y_d) \in \R^d$, we use an affine coordinate system $z
= (\wt z, z_d)$ to represent it so that $z_d= y_d-\Gamma^*(\wt y)$
and $\wt z$ are the coordinates in an orthogonal coordinate system
centered at $x_0$ for the $(d-1)$-dimensional hyperplane $\Pi$ for
the point $(\wt y, \Gamma^*(\wt y))$. Denote such an affine
transformation $y\mapsto z$ by $z=\Psi (y)$. It is clear that there
is a constant $c_5=c_5(\Lambda, R)>1$ so that for every $y\in \R^d$,
$$
c_5 ^{-1}|\wt y -\wt x| \leq |\wt z|\leq c_5 |\wt y -\wt x|, \qquad
c_{5}^{-1} |y-x|\leq |\Psi (y)-\Psi (x)| \leq c_{5} |y-x|
$$
and that
$$
\Psi(E)\subset \{z=(\wt z, z_d)\in \R^d: \ |\widetilde{z}|<c_{5} r_0
\hbox{ and }  0< z_d\leq c_{5} r_0\}.
$$
Denote $x_d-\Gamma^*(\wt x)$ by $w$; that is,
 $\Psi (x)=(\wt 0, w)$.
Hence by \eqref{hhhh} and \eqref{e:KG} and applying the transform
$\Psi$,  we have by using polar coordinates for $\wt z$ on the
hyperplane $\Pi$,
\begin{eqnarray*}
I_3&\leq&c_{6}\int_{ E}\frac{ \overline{h}(y)^{(p-1)_-}\,
|\widetilde{y}-\widetilde{x}|^{2}}{|y-x|^{d+\alpha}}\, dy \leq c_{7}
\int_{\Psi(E)}\ \frac{z_d^
{(p-1)_-} \,
|\wt z|^{2}}{|z- (\wt 0, w)|^{d+\alpha}}\ dz\nonumber\\
&\leq & c_{8}\int_{0}^{c_{5}r_0}  z_d^{(p-1)_-}\left(
\int_0^{c_{5}r_0} \frac{r^{d-2}}{(r+|z_d-w|)^{d+\alpha-2}} dr
\right) dz_d\\
&\leq & c_{8}\int_{0}^{c_{5}r_0}  z_d^{(p-1)_-}\left(
\int_0^{c_{5}r_0}
\frac{1}{(r+|z_d-w|)^{\alpha}} dr \right) dz_d\\
&\leq& c_{9}\int_{0}^{c_{5}r_0} z_d^{(p-1)_-}\left(
\frac1{|z_d-w|^{\alpha-1}} -
\frac{1}{(c_{5}r_0+|z_d-w|)^{\alpha-1}} \right)dz_d\\
&<&  c_{10}\int_{0}^{c_{5}r_0}
\frac1{z_d^{(1-p)^+}\,
|z_d-w|^{\alpha-1}} dz_d  \leq  c_{11}<\infty,
\end{eqnarray*}
where all constants depend on $\alpha$, $p$, $\Lambda$ and $R$. The
last inequality is due to the fact that since $p>0$, $0<\alpha<2$
and $(1-p)^+ + \alpha -1 =\max\{\alpha-p, \alpha -1\}<1$, by
the dominated convergence theorem, $\phi (w):=
\int_0^{c_{5}r_0}
  z_d^{-(1-p)^+}\,  |z_d-w|^{1-\alpha} dz_d$
is a strictly positive continuous function in $x_d\in [0, c_{5}r_0]$
and hence is bounded.  Thus we have proved the claim
\eqref{e:I1234}, hence completing the proof of the lemma. \qed

Since $D$ is a $C^{1,1}$ open set with characteristics $(R,
\Lambda)$, for every
$\lambda \ge 1$, $\lambda D$ is a
$C^{1,1}$ open set
with uniform characteristics $(R, \Lambda)$.
Thus, by the
previous lemma and \eqref{e:L1}, we get the following as a
corollary.

\begin{cor}\label{L:Main1}
Fix $Q \in \partial D$ and the coordinate system $CS_Q$ so that
$$
B(Q, R)\cap D= \big\{ y=(\wt y, \, y_d) \in B(0, R) \mbox{ in }CS_Q:
y_d > \phi (\wt y) \big\}.
$$
Let
$$
h_p(y):=\left(\rho_Q (y)\right)^p{\bf 1}_{D\cap B(Q, 4r_0)}(y) .
$$
Then there exist $C_i=C_i(\alpha, p, \Lambda, R)>0$, $i=5,6,7$,
independent of the choice of the point $Q \in \partial D$ and
$\lambda \ge 1$ such that
\begin{description}
\item {\rm (i)} in the case $\frac{\alpha}2 < p < \alpha $,
for all $x\in D$ such that
$\rho_Q(x)<r_0 \wedge R_*$ and $|\widetilde{x}| <r_0$, we have
 \bee\label{e:h21}
C_6 \left(\rho_Q (x)\right)^{p-\alpha}  \leq \wh \Delta_{d, \lambda}
^{\alpha/2} h_p(x)   \le C_5 \left(\rho_Q (x)\right)^{p-\alpha};
\eee

\item{\rm (ii)}
in the case $ p> {\alpha} $, for
all $x\in D$ such that
$\rho_Q(x)<r_0 \wedge R_*$ and $|\widetilde{x}| <r_0$, we have
\bee\label{e:h01} |\wh \Delta_{d, \lambda}^{\alpha/2} h_p(x)|  \le
C_7 \lambda^{p-\alpha}; \eee

\item{\rm (iii)}
in the case $ p= {\alpha} $, for
all $x\in D$ such that
$\rho_Q(x)<r_0 \wedge R_*$ and $|\widetilde{x}| <r_0$, we have
\bee\label{e:h01_1} |\wh \Delta_{d, \lambda}^{\alpha/2} h_p(x)|  \le
C_7  \left|\log \left(\rho_Q (x)/\lambda\right) \right|\, .
\eee

\end{description}
\end{cor}

The following scaling property of $X^{a}$ will be used below:
If $(X^{a, D}_t, t\geq 0)$ is the  subprocess in $D$ of the
independent sum of a Brownian motion and a symmetric $\alpha$-stable
process on $\R^d$ with weight $a$, then $(\lambda X^{a,
D}_{\lambda^{-2} t}, t\geq 0)$ is the subprocess in $\lambda D$ of
the independent sum of a Brownian motion and a symmetric
$\alpha$-stable process on $\R^d$  with  weight $a
\lambda^{(\alpha-2)/\alpha}$. So for any $\lambda>0$, we have
\begin{equation}\label{e:scaling}
p^{a\lambda^{(\alpha-2)/\alpha}}_{\lambda D} ( t,  x, y)=
\lambda^{-d} p^{a}_D (\lambda^{-2}t, \lambda^{-1} x, \lambda^{-1} y) \qquad
\hbox{for } t>0 \hbox{ and } x, y \in
\lambda D.
\end{equation}
By integrating the above equation with respect to $t$, we get
\begin{equation}\label{e:scaling2}
G^{a\lambda^{(\alpha-2)/\alpha}}_{\lambda D} (x, y)=
\lambda^{2-d} G^{a}_D ( \lambda^{-1} x, \lambda^{-1} y) \qquad \hbox{for }
x, y \in \lambda D
\end{equation}
where
$$  G^a_D(x, y):=\int_0^\infty
p^a_D(t, x, y)dt
$$
is the Green function of $X^{a}$ in $D$. It is well known that the
L\'evy measure of $X^1$ has the intensity
$$
J^1(x, y)=j^1(|x-y|)={\cal A}(d, \alpha)|x-y|^{-(d+\alpha)}.
$$
Thus by a scaling argument, we get that the L\'evy intensity of
$X^a$ is
$$
J^a(x, y)=j^a(|x-y|)=a^{\alpha}{\cal A}(d,
\alpha)|x-y|^{-(d+\alpha)},
$$
which gives the L\'evy system \eqref{e:levy} of $X^a$.

By a $\lambda$-truncated symmetric $\alpha$-stable process in $\R^d$
we mean a pure jump symmetric L\'evy process
$\wh Y^\lambda=(\wh Y^\lambda_t, t \ge 0, \P_x )$
in $\R^d$  with   L\'evy density ${\cal A}(d,
\alpha) |x|^{-d-\alpha}\, \1_{\{|x|< \lambda\}}$.  Note that the
L\'evy exponent $\psi^\lambda$ of $\wh Y^\lambda$ defined by
$$
\E_x\left[e^{i\xi\cdot(\wh Y^\lambda_t-\wh Y^\lambda_0)}\right] =
e^{-t\psi^\lambda(\xi)} \quad \quad \mbox{ for every } x\in \R^d
\mbox{ and } \xi\in \R^d,
$$
is given  by
\begin{equation}\label{e:psi}
\psi^\lambda(\xi)={\cal A}(d, \alpha) \int_{\{|y|<
\lambda\}}\frac{1-\cos(\xi\cdot y)}{|y|^{d+\alpha}}dy.
\end{equation}

Suppose that $\wh Y^{\lambda/a}$ is a $(\lambda/a)$-truncated
symmetric $\alpha$-stable process in $\R^d$ which is independent of
the Brownian motion $X^0$. For any $a>0$, we define
$$
\wh X_t^{a,\lambda}:=X^0_t+ a\wh Y^{\lambda/a}_t, \qquad t\geq 0 .
$$
Note that from
\eqref{e:psi} we can easily check that for any $b>0$,
\begin{equation}\label{e:ps}
\psi^\lambda(b\xi)=b^{\alpha}\psi^{\lambda b}(\xi)\, \qquad  \hbox{ for every }\xi\in \R^d.
\end{equation}
Thus
for any $a>0$ and $\xi, x \in \R^d$,
$$
\E_x\left[e^{i\xi\cdot(\wh X_t^{a,\lambda}-\wh
X_0^{a,\lambda})}\right] = e^{-t|\xi|^2}
\E_x\left[e^{i(a\xi)\cdot(\wh Y^{\lambda/a}_t - \wh
Y^{\lambda/a}_0)}\right]= e^{-t(|\xi|^2+\psi^{\lambda/a}(a\xi))}=
e^{-t(|\xi|^2+a^\alpha\psi^{\lambda}(\xi))}.
$$
Therefore $\wh X^{a, \lambda}$ has the same distribution as the
L\'evy process obtained from $X^a$ by removing jumps of size larger
than $\lambda$. The above observation also gives us that the
infinitesimal generator of $\wh X^{a, \lambda}$ is $\Delta +
a^\alpha \wh \Delta_{d, \lambda}^{\alpha/2}$, and the L\'evy
intensity for $\wh X^{a, \lambda}$ is
$$
J^{a, \lambda}(x, y):=a^{\alpha}{\cal A}(d, \alpha)
|x-y|^{-(d+\alpha)} \, \1_{\{|x-y| < \lambda \}}\, .
$$
The L\'evy intensity describes the jumps of the process $\wh X^{a, \lambda}$
through the L\'evy system:
for any non-negative measurable
function $f$ on $\bR_+ \times \bR^d\times \bR^d$, $x\in \bR^d$ and
stopping time $T$ (with respect to the filtration of $\wh
X^{a,\lambda}$),
\begin{equation}\label{e:levyt}
\E_x \left[\sum_{s\le T} f(s,\wh X^{a, \lambda}_{s-}, \wh X^{a,
\lambda}_s) \right]= \E_x \left[ \int_0^T \left( \int_{\bR^d}
f(s,\wh X^{a, \lambda}_s, y) J^{a, \lambda}(\wh X^{a, \lambda}_s,y)
dy \right) ds \right].
\end{equation}
For any open set $U\subset \bR^d$, let $\wh
\tau^{a,\lambda}_U=\inf\{t>0: \, \wh X^{a,\lambda}_t\notin U\}$ be
the first exit time from $U$ by $\wh X^{a,\lambda}$, and denote by
$\wh X^{a, \lambda, U}$ the subprocess of $\wh X^{a, \lambda}$
killed upon leaving $U$.
When $\lambda=1$, we simply write $\wh
X^{a}$ for $\wh X^{a,1}$ and $\wh \tau^{a}_U$ for $\wh
\tau^{a,1}_U$. The following scaling property will be used in the
next lemma: by \eqref{e:ps}, we see that for every $\lambda, a, b>0$
and $\xi, x \in \R^d$,
\begin{align*}
&\E_x\left[e^{i\xi\cdot(b (\wh X_{b^{-2} t}^{a,\lambda}-\wh
X_0^{a,\lambda}))}\right] = e^{-t|\xi|^2} \E_x\left[e^{i(ab\xi)
\cdot(\wh Y^{\lambda/a}_{b^{-2} t}-\wh Y^{\lambda/a}_0)}\right]\\
&=e^{-t|\xi|^2} e^{-b^{-2} t \psi^{\lambda/a}(ab
\xi)}=e^{-t(|\xi|^2+a^\alpha b^{\alpha-2}\psi^{b\lambda}(\xi))}.
\end{align*}
Thus, if $\{\wh X^{a, \lambda, D}_t, t\geq 0\}$ is the subprocess of
$\{\wh X^{a, \lambda}_t, t\geq 0\}$ in $D$, then $ \{b \wh X^{a,
\lambda, D}_{b^{-2} t}, t\geq 0\}$ is the subprocess of $\{\wh
X^{ab^{(\alpha-2)/\alpha}, b\lambda}_t, t\geq 0\}$  in $b D$. In
particular, if $\{\wh X^{a, D}_t, t\geq 0\}$ is the subprocess  of
$\{\wh X^{a}_t, t\geq 0\}$ in $D$, then $ \{\lambda \wh X^{a,
D}_{\lambda^{-2} t}, t\geq 0\}$ is the subprocess of $\{\wh
X^{a\lambda^{(\alpha-2)/\alpha}, \lambda}_t, t\geq 0\}$ in $\lambda
D$.
So for any $\lambda>0$, we have
\begin{equation}\label{e:scaling_a}
\wh p^{\ a\lambda^{(\alpha-2)/\alpha}, \lambda}_{\lambda D} (
t, x, y)= \lambda^{-d} \wh p^{\ a,1}_D (\lambda^{-2}t, \lambda^{-1} x,
\lambda^{-1} y) \qquad \hbox{for } t>0 \hbox{ and } x, y \in \lambda
D
\end{equation}
where $\wh p^{\ a, \lambda}_{D} ( t,  x, y)$ is the transition density
of $\wh X^{a, \lambda,D}$. By integrating the above equation with
respect to $t$, we get
\begin{equation}\label{e:scaling2_a}
\wh G^{\, a\lambda^{(\alpha-2)/\alpha},\lambda}_{\lambda D}
(x, y)= \lambda^{2-d} \wh G^{\, a,1}_D ( \lambda^{-1} x, \lambda^{-1} y)
\qquad \hbox{for } x, y \in \lambda  D
\end{equation}
where
$$
\wh  G^{\, a, \lambda}_D(x, y):=\int_0^\infty \wh p^{\ a,\lambda}_D(t, x,
y)dt
$$
is the Green function of $\wh X^{a,\lambda}$ in $D$.

\medskip

For our reader's convenience, we summarize some notations below.

\halign{{}\hskip 0.6truein  #\hfil \hskip 0.6truein   &#\hfil \hskip
0.6truein &#\hfill \cr {\it Process} & {\it Generator}  & {\it
L\'evy (jumping) kernel} \medskip \cr $X^0$ &  $\Delta$ & 0 \cr $Y$
& $\Delta^{\alpha/2}$  & $\sA(d, \alpha )\, |z|^{-d-\alpha}$
\medskip\cr $aY$ &  $a^\alpha \Delta^{\alpha/2}$  & $a^\alpha \sA(d,
\alpha )\, ||z|^{-d-\alpha}$
   \medskip \cr
$\wh Y^\lambda$ &  $\wh \Delta_{d, \lambda}^{\alpha/2}$
   & $\sA(d, \alpha)\, |z|^{-d-\alpha}\, \1_{\{|z|<\lambda\}}$ \medskip \cr
 $X^a:=X^0+aY$ &  $\Delta+a^\alpha \Delta^{\alpha/2}$ &
 $ a^\alpha \sA(d, \alpha)\, |z|^{-d-\alpha} $  \medskip\cr
$\wh X^{a, \lambda}:=X^0+a \wh Y^{\lambda/a}$& $\Delta+a^\alpha \wh
\Delta_{d, \lambda}^{\alpha/2}$ & $a^\alpha \sA(d, \alpha)\,
|z|^{-d-\alpha}\, \1_{\{|z|<\lambda\}}$  \medskip\cr $\wh X^{a
}:=\wh X^{a, 1}$& $\Delta+a^\alpha \wh \Delta_{d }^{\alpha/2}$ &
$a^\alpha \sA(d, \alpha)\, |z|^{-d-\alpha}\, \1_{\{|z|<1\}}$ . \cr }

\bigskip

Recall that $\rho_Q (x) := x_d -  \phi_Q (\wt x) $ for every $Q\in
\partial D$ and $ x \in \left\{ y=(\wt y, y_d)\in B(Q, R): y_d > \phi_Q
(\wt y) \right\}$. We define for $r_1, r_2>0$
$$
D_Q( r_1, r_2) :=\left\{ y\in D: r_1 >\rho_Q(y) >0,\, |\wt y | < r_2
\right\}.
$$

\begin{lemma}\label{L:2}
 There are constants
$\delta_0=\delta_0( R, M, \Lambda, \alpha)\in (0, r_0)$, $C_8=C_8(
R, M, \Lambda, \alpha)>0$ and $C_9=C_9( R, M, \Lambda, \alpha)>0$
such that for every $a \in (0, M]$,
$\lambda \ge 1$,
$Q \in
\partial D$ and $x \in D_Q (\lambda^{-1} \delta_0, \lambda^{-1} r_0)$ with
$\wt x =0$,
 \bee\label{e:L:1}
\P_{x}\left(\wh X^{a}_{\wh \tau^{a}_{D_Q (\lambda^{-1} \delta_0 ,
\lambda^{-1} r_0 )}} \in  D_Q (2\lambda^{-1} \delta_0 ,\lambda^{-1} r_0 ) \right)
\ge\,{C_8}{\lambda} \delta_D (x),
 \eee
 \bee\label{e:L:2}
\P_{x}\left(\wh X^{a}_{\wh \tau^{a}_{ D_Q (\lambda^{-1} \delta_0 ,
\lambda^{-1} r_0 )}} \in  D\right)\le {C_9}{\lambda} \delta_D (x)
 \eee
and
 \bee\label{e:L:3}
\E_x\left[\wh \tau^{a}_{ D_Q (\lambda^{-1} \delta_0 , \lambda^{-1}
r_0)}\right]\,\le\, {C_9}{\lambda^{-1}} \delta_D (x).
 \eee
\end{lemma}

\pf
To derive the estimates in the lemma, it will be convenient to
consider the scaled process $\lambda \wh X^a_{\lambda^{-2}t}$, which
has the same distribution as $ \wh X^{a\lambda^{(\alpha -2)/\alpha},
\lambda}$. The latter has infinitesimal generator $\Delta + a^\alpha
\lambda^{\alpha -2} \wt \Delta_{d, \lambda}^{\alpha/2}$.

Without loss of generality, we assume $Q=0$ and let $\phi:
\bR^{d-1}\to \bR$ be the $C^{1,1}$-function satisfying $\phi (\wt
0)= \nabla\phi (\wt 0)=0$, $\| \nabla \phi  \|_\infty \leq \Lambda$,
$| \nabla \phi (\wt y)-\nabla \phi (\wt z )| \leq \Lambda |\wt y-
\wt z|$ and $CS_Q$ be the corresponding coordinate system  such that
$$
B(Q, R)\cap D= \big\{(\wt y, \, y_d) \in B(0, R)\textrm{ in } CS_Q: y_d >
\phi (\wt y) \big\}.
$$
Note that, since $D$ is a $C^{1,1}$ open set with characteristics
$(R, \Lambda)$,
for every $\lambda\ge 1$, $\lambda D$ is a
$C^{1,1}$ open set with
 the same characteristics $(R, \Lambda)$.
Let
$\phi_{\lambda }(\wt y) :=\phi(\lambda^{-1} \wt y): \bR^{d-1}\to \bR$.
Then $\phi_{\lambda}$ satisfies $\phi_{\lambda} (\wt 0)=
\nabla\phi_{\lambda} (\wt 0)=0$, $\| \nabla \phi_{\lambda} \|_\infty
\leq  \Lambda$, $| \nabla \phi_{\lambda} (\wt y)-\nabla
\phi_{\lambda} (\wt z)| \leq  \Lambda |\wt y- \wt z|$ and
$$
B(Q, R)\cap {\lambda}D=  \big\{ y\in B(0, R) \mbox{ in } CS_Q :
\ y_d > \phi_{\lambda} (\wt y) \big\} \quad \text{ for all }
\lambda \ge 1\, .
$$
We let $p>0$ be such that $p\not=\alpha$ and
$1 < p< \left( 2 \wedge (3-\alpha) \right) $,
 and define
\begin{eqnarray*}
\rho_{\lambda} (y) &:=& y_d -  \phi_{\lambda} (\wt y), \\
h_{\lambda}(y) &:=& \rho_{\lambda}(y)
{\bf 1}_{B(0, 4r_0) \cap \lambda  D}(y), \\
h_{\lambda, p}(y) &:=& h_\lambda (y)^p \,=\,
\left(\rho_{\lambda}(y)\right)^p
{\bf 1}_{B(0, 4r_0) \cap \lambda  D}(y), \\
D (\lambda, r_1, r_2) &:=& \left\{ y\in \lambda  D: \,
0<\rho_{\lambda}(y) <r_1 \hbox{ and  }  |\wt y | < r_2 \right\}.
\end{eqnarray*}
Since $\rho_\lambda (y) \leq \sqrt{1+\Lambda^2} \,
\delta_{\lambda  D} (y)$ in view of \eqref{e:d_com}, we have
$0\leq h_\lambda \leq R\leq 1$. It is easy to see that $D (\lambda,
r_1, r_2)$ is contained in $D \cap B(0, R/4)$ for every $r_1, r_2
\le r_0$. Note that the (vector-valued) Lipschitz function $\nabla
\phi_\lambda$ is differentiable almost everywhere. So for a.e.
$y \in B(0, 4r_0)\cap \lambda  D$,
\begin{align}\label{e:**}
\Delta h_{\lambda}(y) = \Delta (y_d -  \phi_{\lambda} (\wt y))
=-\Delta \phi_{\lambda}(\wt y)
\end{align}
and
\begin{eqnarray*}
\Delta h_{\lambda, p}(y) &=& \Delta (y_d -  \phi_{\lambda} (\wt y))^p \\
&=&p(p-1)(1+ |\nabla \phi_{\lambda} (\wt y)|^2)(\rho_{\lambda}
(y))^{p-2} - p\, (\rho_{\lambda} (y))^{p-1} \Delta \phi_{\lambda}(\wt y)\\
&\ge & p(p-1)(1 + |\nabla \phi_{\lambda} (\wt y)|^2) (\rho_{\lambda}
(y))^{p-2}-p\, (\rho_{\lambda} (y))^{p-1} \| \Delta
\phi_{\lambda}\|_\infty.
\end{eqnarray*}
Thus, since $p\in (1,2)$,  we can choose a positive constant
$\delta_1=\delta_1(R,M,\Lambda,\alpha)\in (0,r_0)$,
independent of $\lambda$, so that there is $c_1>0$ such that
\bee\label{e:delta} \Delta h_{\lambda, p}(y) \ge c_1 (\rho_{\lambda}
(y))^{p-2} >0 \qquad \hbox{for a.e. } y \in
D(\lambda,\delta_1,r_0)\, . \eee

We divide the rest of the proof into three steps.

\medskip
\noindent{\it Step 1: Constructing suitable superharmonic and
subharmonic functions with respect to $\Delta + a^\alpha
\lambda^{\alpha -2} \wt \Delta_{d, \lambda}^{\alpha/2}$}. Let $\psi$
be a smooth positive function on $\bR^d$ with bounded first and
second order partial derivatives such that $\psi(y)= 2^{p+1}
|\widetilde y|^2/ r_0^2$ for $|y| <r_0/4$ and $2^{p+1} \leq \psi (y)
\leq 2^{p+2}$ for $|y| \ge r_0/2$. Now we consider
$$
u_{1, \lambda}(y):=h_{\lambda}(y) +  h_{\lambda, p}(y)
$$
and
$$u_{2, \lambda}(y):=h_{\lambda}(y) +  \psi(y) -h_{\lambda, p}(y).
$$
Observe that since $0\leq h_{\lambda}\leq 1$ and $p\geq 1$,
both $u_{1, \lambda}$ and $u_{2, \lambda}$ are non-negative.
By Taylor's expansion with remainder of order $2$,
  \bee\label{e:tr2}
\left| \left(\Delta+ a^\alpha \lambda^{(\alpha-2)}\wh \Delta_{d,
\lambda}^{\alpha/2}\right)\psi(y)\right| \le|\Delta \psi(y) |
+M^\alpha \left|\wh \Delta_{d, \lambda}^{\alpha/2}\psi(y)\right| \le
c_2(\alpha, M)< \infty.
 \eee
Note that the constant $c_2$ above is independent of $\lambda$.
Moreover,
 since $\lambda \geq 1$, $p>\alpha/2$ and $p\not= \alpha$,
by \eqref{e:h21} and \eqref{e:h01} there exist
$c_3=c_3( R, \Lambda)>0$ and $\delta_2=\delta_2( R, \Lambda) \in (0, \delta_1]$
independent of $\lambda$ such that
$$
\wh \Delta_{d, \lambda}^{\alpha/2}h_{\lambda, p}(y) \ge - c_3
\lambda^{p-\alpha}
\quad \text{for } y \in D(\lambda, \delta_2,r_0).
$$
Thus by using \eqref{e:delta}, the fact that
$p<2 $
 and the inequality above, and by choosing
$\delta_2$ smaller if necessary, we get
  \begin{align}\label{e:gp}
\left(\Delta+ a^\alpha \lambda^{(\alpha-2)} \wh \Delta_{d,
\lambda}^{\alpha/2} \right)h_{\lambda, p}(y) &\ge
c_1\rho_{\lambda}(y)^{p-2} -M^{\alpha} c_3 \lambda^{(p-\alpha)+
(\alpha-2)}\\
 &\ge c_1 \rho_{\lambda}(y)^{p-2} -M^{\alpha} c_3 \ge
\frac{c_1}2 \rho_{\lambda}(y)^{p-2} \nonumber
 \end{align}
for a.e.~$ y \in D(\lambda, \delta_2,r_0).$ Furthermore by
 \eqref{e:h21}-\eqref{e:h01_1}  and \eqref{e:**},
 there exist $c_4=c_4(
M)>0$ and $\delta_3 \in (0, \delta_2)$ independent of
$\lambda \ge 1$
such
that for a.e.~$y \in D(\lambda, \delta_3,r_0)$
  \begin{eqnarray}\label{e:gp2}
&&\left| \left(\Delta+ a^\alpha \lambda^{(\alpha-2)} \wh \Delta_{d,
\lambda}^{\alpha/2} \right) h_{\lambda}(y) \right| \nonumber\\
 &\le& c_4
\left(1+\lambda^{(1-\alpha)^++
(\alpha-2)}\rho_{\lambda}(y)^{(1-\alpha)\wedge 0}
+ \1_{\{\alpha =1\}}\, \lambda^{-1} |\log (\rho_\lambda(y)/\lambda)|  \right) \\
 &\le& c_4
\left(1+\lambda^{(1-\alpha)^++
(\alpha-2)}\rho_{\lambda}(y)^{(1-\alpha)\wedge 0}
 + e^{-1}+\1_{\{\alpha =1\}}\,   |\log  \rho_\lambda(y) |  \right). \nonumber
 \end{eqnarray}
Thus by  \eqref{e:tr2}-\eqref{e:gp2} and the fact that $p<2 \wedge
(3-\alpha)$, there exists $\delta_4 \in (0, \delta_3)$ independent
of $\lambda \ge 1$ such that
 \bee\label{e:sh}
\left(\Delta+  a^\alpha \lambda^{(\alpha-2)} \wh \Delta_{d,
\lambda}^{\alpha/2}\right)u_{2, \lambda}(y) \le c_2 +
c_4 \left(2+  |\log \rho_\lambda(y)|  +
\rho_{\lambda}(y)^{(1-\alpha)\wedge 0}\right) - \frac{c_1}{2}
\rho_{\lambda}(y)^{p-2} \le -1
 \eee
for a.e.~$ y \in D(\lambda, \delta_4,r_0)$.

On the other hand, we have from \eqref{e:h01} and \eqref{e:h01_1},
  \begin{eqnarray*}
  \left(\Delta+ a^\alpha \lambda^{(\alpha-2)} \wh \Delta_{d,
\lambda}^{\alpha/2}\right)h_{\lambda}(y)
&\ge&
  -\|\Delta \phi_\lambda
\|_{\infty} -c_5 M^{\alpha}
(\lambda^{(1-\alpha) + (\alpha-2)} +\lambda^{-1} \log \lambda    +\lambda^{-1} |\log \rho_\lambda(y)|) \\
  &\ge&-\|\Delta \phi_\lambda \|_{\infty}
-c_5 M^{\alpha} (1+ e^{-1} +|\log \rho_\lambda(y)| )
 \end{eqnarray*}
for a.e.~$y\in D(\lambda, \delta_4,r_0).$ Combining the inequality
above with \eqref{e:gp}, by choosing $\delta_4$ smaller if
necessary, we have for a.e.~$ y \in D(\lambda, \delta_4,r_0),$
 \bee\label{e:sh2}
\left(\Delta+ a^\alpha \lambda^{(\alpha-2)} \wh \Delta_{d,
\lambda}^{\alpha/2}\right)u_{1, \lambda}(y) \ge -\|\Delta \phi
\|_{\infty}
  -c_5 M^{\alpha}   (2+  |\log \rho_\lambda(y)|)
+ \frac{c_1}2\rho_{\lambda}(y)^{p-2} \ge
0.
 \eee

\medskip
\noindent{\it Step 2: Translating super-/sub-harmonic functions into
super-/sub-martingale properties for  $ \wh
X^{a\lambda^{(\alpha-2)/\alpha},\, \lambda}$}. For notational
convenience, we let
$$
\wt X^{a,\lambda} :=\wh X^{a\lambda^{(\alpha-2)/\alpha},
\lambda}\quad \text{and}\quad \wt \tau_U^{\, a,\lambda}:=
\wh\tau_U^{\, a\lambda^{(\alpha-2)/\alpha}, \lambda}.
$$
We claim that the estimates \eqref{e:sh} and \eqref{e:sh2} imply
that
\begin{equation}\label{e:claim1}
t\mapsto u_{2, \lambda} \Big( \wt X^{a,\lambda}_{t\wedge \wt
\tau^{a, \lambda}_{D(\lambda, \delta_4,r_0)}}\Big) \  \hbox{ is a
bounded supermartingale},
\end{equation}
\begin{equation}\label{e:claim3}
\E_x \left[ \wt \tau^{a, \lambda}_{D(\lambda, \delta_4,r_0)} \right]
\leq \rho_{\lambda} (x),
\end{equation}
and
\begin{equation}\label{e:claim2}
t\mapsto u_{1, \lambda} \Big( \wt X^{a,\lambda}_{t\wedge \wt
\tau^{a, \lambda}_{D(\lambda, \delta_4,r_0)}}\Big) \ \hbox{ is a
bounded submartingale}.
\end{equation}

Observe that if $v$ is a bounded $C^2$-function on $\R^d$ with
bounded second order partial derivatives, then by Ito's formula and
the L\'evy system \eqref{e:levyt},
\begin{equation}\label{e:cl0}
M^v_t=v(\wt X^{a,\lambda}_t)-v(\wt X^{a,\lambda}_0)-\int_0^t
\left(\Delta+ a^\alpha \lambda^{(\alpha-2)} \wh \Delta_{d,
\,\lambda}^{\alpha/2}\right)v(\wt X^{a,\lambda}_s)ds
\end{equation}
is a martingale (see the proof of
Proposition 4.1 in \cite{BC} for the derivation of a similar
assertion). If the functions $u_{2, \lambda}$ and $u_{1, \lambda}$
were $C^2$ with bounded second order partial derivatives, then the
claims \eqref{e:claim1}, \eqref{e:claim3} and \eqref{e:claim2} would
just follow from \eqref{e:cl0} and the estimates \eqref{e:sh} and
\eqref{e:sh2}. However they are not $C^2$ since $D$ is $C^{1,1}$ and
they are truncated on the outside of $B(0, 4r_0) \cap \lambda D$.
So we will use a mollifier.
Let $g$ be a non-negative smooth
function with compact support in $\R^d$ whose value only depends on
$|x|$ such that $g(x)=0$ for $|x|>1$ and $\int_{\R^d} g (x)
dx=1$. For $k\geq 1$, define $g_k(x)= 2^{kd} g (2^k x)$. Set
$$
u_{i, \lambda}^{(k)}(z):= ( g_k*u_{i, \lambda})(z) :=\int_{\R^d}
g_k (y) u_{i, \lambda}(z-y)dy, \qquad i=1,2.
$$
As
$$
\left(\Delta+ a^\alpha \lambda^{(\alpha-2)} \wh \Delta_{d,
\,\lambda}^{\alpha/2}\right)u^{(k)}_{i, \lambda} = g_k*
\left(\Delta+ a^\alpha \lambda^{(\alpha-2)} \wh \Delta_{d,
\,\lambda}^{\alpha/2}\right) u_{i, \lambda} \qquad \hbox{for } i=1,
2,
$$
we have by \eqref{e:sh} and \eqref{e:sh2} that
$$
\left(\Delta+ a^\alpha \lambda^{(\alpha-2)} \wh \Delta_{d, \,\lambda
}^{\alpha/2}\right)u^{(k)}_{1, \lambda}\geq 0 \quad \hbox{and} \quad
\left(\Delta+ a^\alpha \lambda^{(\alpha-2)} \wh \Delta_{d, \,\lambda
}^{\alpha/2}\right)u^{(k)}_{2, \lambda}\leq -1
$$
on $D_k(\lambda, \delta_4, r_0):=\left\{y: \delta_4 -2^{-k}> \rho_\lambda(y)
>2^{-k} \hbox{ and } |\wt y|<r_0-2^{-k}\right\}$.

Since $u^{(k)}_{i,
\lambda}$, $i=1, 2$, are bounded smooth functions on $\R^d$ with
bounded first and second order partial derivatives, it follows from
\eqref{e:cl0} that
$$
t\mapsto u^{(k)}_{2, \lambda} \Big( \wt X^{a,\lambda}_{t\wedge \wt
\tau^{a, \lambda}_{D_k(\lambda, \delta_4,r_0)}}\Big)
 + t\wedge \wt \tau^{a, \lambda}_{D_k(\lambda, \delta_4,r_0)}
 \ \hbox{ is a positive supermartingale}
$$
and
$$
t\mapsto u^{(k)}_{1, \lambda} \Big( \wt X^{a,\lambda}_{t\wedge \wt
\tau^{a, \lambda}_{D_k(\lambda, \delta_4,r_0)}}\Big) \ \hbox{ is a
bounded submartingale}.
$$
Since for $i=1, 2$, $u_{i, \lambda}$ is bounded and continuous,
$u^{(k)}_{i, \lambda}$ converges uniformly to  $u_{i, \lambda}$.
Thus
\begin{equation}\label{e:3.37}
t\mapsto u_{2, \lambda} \Big( \wt X^{a,\lambda}_{t\wedge \wt
\tau^{a, \lambda}_{D_k(\lambda, \delta_4,r_0)}}\Big) + t\wedge \wt
\tau^{a, \lambda}_{D_k(\lambda, \delta_4,r_0)} \ \hbox{is a positive
supermartingale  }
\end{equation}
and
$$
t\mapsto u_{1, \lambda} \Big( \wt X^{a,\lambda}_{t\wedge \wt
\tau^{a, \lambda}_{D_k(\lambda, \delta_4,r_0)}}\Big) \ \hbox{ is a
bounded submartingale}.
$$
Since $D_k(\lambda, \delta_4,r_0)$ increases to $D(\lambda,
\delta_4,r_0)$, we conclude that \eqref{e:claim1} and
\eqref{e:claim2} hold. Moreover,  for each fixed $k\geq 1$ and
$t>0$, we have from \eqref{e:3.37} that
$$
\E_x\left[  u_{2, \lambda} \Big( \wt X^{a,\lambda}_{t\wedge \wt
\tau^{a, \lambda}_{D_k(\lambda, \delta_4,r_0)}}\Big) + t\wedge \wt
\tau^{a, \lambda}_{D_k(\lambda, \delta_4,r_0)}\right] \leq u_{2,
\lambda}(x).
$$
Since $u_{2, \lambda}\geq 0$, by first letting $k\to \infty$ and
then $t\to \infty$, we get $\E_x\left[  \wt \tau^{a, \lambda}_{D
(\lambda, \delta_4,r_0)}\right] \leq u_{2, \lambda}(x)$. Since $\wt
x=0$, $\psi (x)=0$ and so $u_{2, \lambda}(x)\leq \rho_\lambda(x)$.
This proves \eqref{e:claim3}.

\medskip
\noindent{\it Step 3: Deriving the desired exit distribution
estimates by utilizing the super-/sub-martingale property}. Since
$\psi \ge 2^{p+1}$ on $|\wt y| \ge r_0$ and $\psi(x)=0$, we have by
\eqref{e:claim1},
\begin{align*}
\rho_{\lambda}(x) &\ge u_{2, \lambda}(x)\\
&\ge \E_{x}\left[u_{2, \lambda}\left(\wt X^{a,\lambda}_{\wt \tau^{a,
\lambda}_{D(\lambda, \delta_4,r_0)}} \right); \wt X^{a,\lambda}_{\wt
\tau^{a, \lambda}_{D(\lambda,\delta_4,r_0)}} \in
(\lambda  D) \setminus  D(\lambda, \infty,r_0)  \right]\\
&\ge
(2^{p+1}-1) \,
 \P_{x}\left(\wt X^{a,\lambda}_{\wt \tau^{a,
\lambda}_{D( \lambda, \delta_4,r_0)}} \in (\lambda  D) \setminus
D( \lambda, \infty,r_0) \right) .
\end{align*}
We also have from \eqref{e:claim2}
\begin{align*}
\rho_{\lambda}(x) &\le \rho_{\lambda}(x)+\rho_{\lambda}(x)^{p}=
u_{1, \lambda}(x) \le \E_{x}\left[u_{1, \lambda} \left(\wt X^{a,
\lambda}_{ \wt \tau^{a, \lambda}_{D(\lambda, \delta_4,r_0)}} \right)
\right] \le 2 \P_{x}\left(\wt X^{a,\lambda}_{\wt \tau^{a,
\lambda}_{D(\lambda,\delta_4,r_0)}} \in \lambda  D \right) .
\end{align*}
Combining the two displays above, we get
\begin{align}
&\P_{x}\left(\wt X^{a,\lambda}_{\wt \tau^{a, \lambda}_{D(\lambda,
\delta_4,r_0)}} \in  D(\lambda, \infty,r_0)
\right)\label{e:ss}\\
=&\P_{x}\left(\wt X^{a,\lambda}_{\wt \tau^{a,
\lambda}_{D(\lambda,\delta_4,r_0)}} \in \lambda  D  \right)-
\P_{x}\left(\wt X^{a,\lambda}_{\wt \tau^{a, \lambda}_{D(\lambda,
\delta_4,r_0)}} \in (\lambda  D) \setminus  D(\lambda,
\infty,r_0)
\right)\nonumber\\
\geq & \frac{2^{p+1}-3}{2(2^{p+1}-1)}
\rho_{\lambda}(x). \nonumber
\end{align}
By \eqref{e:levyt},
\begin{align}
&\P_{x}\left(\wt X^{a,\lambda}_{\wt \tau^{a, \lambda}_{D(\lambda,
\delta_4,r_0)}} \in  D(\lambda, \infty,r_0) \setminus D(\lambda,
2\delta_4 ,r_0) \right)\label{e:ss2}\\
=&\E_x \left[ \int_0^{\wt \tau^{a, \lambda}_{D(\lambda,
\delta_4,r_0)}} \int_{  D(\lambda, \infty,r_0) \setminus D(\lambda,
2\delta_4 ,r_0)}\frac{(a\lambda^{(\alpha-2)/\alpha})^{\alpha}\sA (d,
\alpha) }{|\wt X^{a,\lambda}_s-y|^{d+\alpha} } \1_{\{|\wt
X^{a,\lambda}_s-y|< \lambda\}}dy ds\right]\nonumber\\
\le&\E_x \left[ \int_0^{\wt \tau^{a, \lambda}_{D(\lambda,
\delta_4,r_0)}} \int_{  D(\lambda, \infty,r_0) \setminus D(\lambda,
2\delta_4 ,r_0) }\frac{a^{\alpha}\lambda^{\alpha-2}\sA (d,
\alpha)}{|\wt X^{a,\lambda}_s-y|^{d+\alpha} }dy ds\right]
\nonumber\\
\le& c_6\sA (d, \alpha) a^{\alpha}\lambda^{\alpha-2}\left(\int_{
D(\lambda, \infty, r_0) \setminus D(\lambda, 2\delta_4, r_0)}
|y|^{-d-\alpha} dy \right)\,\E_x \left[\wt \tau^{a,
\lambda}_{D(\lambda, \delta_4,r_0)}
\right]  \nonumber\\
\le& c_7\sA (d, \alpha) a^{\alpha}\lambda^{\alpha-2}
\left(\int_{D(\lambda, 2\delta_4 ,r_0) \setminus D(\lambda,
3\delta_4/2 ,r_0)}|y|^{-d-\alpha} dy \right) \, \E_x \left[\wt
\tau^{a, \lambda}_{D(\lambda, \delta_4,r_0)}\right]
\nonumber\\
\leq & c_8 \, \E_x \left[ \int_0^{\wt \tau^{a, \lambda}_{D(\lambda,
\delta_4,r_0)}} \int_{ D(\lambda, 2 \delta_4,r_0) \setminus
D(\lambda, 3\delta_4/2 ,r_0)}\frac{a^{\alpha}\lambda^{\alpha-2}\sA
(d, \alpha)} {|\wt X^{a,\lambda}_s-y|^{d+\alpha} } \, \1_{\{|\wt
X^{a,\lambda}_s-y|< \lambda\}}dy ds
\right]\nonumber\\
=& c_8 \, \P_{x}\left(\wt X^{a,\lambda}_{\wt \tau^{a,
\lambda}_{D(\lambda, \delta_4, r_0)}} \in D(\lambda,2\delta_4,r_0)
\setminus D( \lambda, 3\delta_4/2 ,r_0) \right). \nonumber
\end{align}
Thus from \eqref{e:ss}-\eqref{e:ss2}
 \bee\label{e:lbp}
\P_{x}\left(\wt X^{a,\lambda}_{\wt \tau^{a, \lambda}_{D(\lambda,
\delta_4,r_0)}} \in D(\lambda,2\delta_4,r_0) \right)\ge c_9
\rho_{\lambda}(x).
 \eee

Recall that $0\leq h_{\lambda, p}\leq 1$. If $|y|>r_0/2$, then
$\psi(y)\ge 2^{p+1}$, we have
$$
u_{2,\lambda}(y)=\psi(y)+h_{\lambda}(y)-h_{\lambda,p}(y)\ge \psi(y)-
h_{\lambda,p}(y)\ge 2^p \ge 1  \qquad \hbox{for } y\in B(0,r_0/2)^c\, .
$$
Furthermore, for $y\in B(0,4r_0)$ such that $\delta_4\le \rho_{\lambda}(y)<4r_0$,
$$
u_{2,\lambda}(y)=\psi(y)+h_{\lambda}(y)-h_{\lambda,p}(y)\ge
\rho_{\lambda}(y)-\rho_{\lambda}(y)^p \ge c_{10}\, ,
$$
where $c_{10}\in (0,1)$ depends on $\delta_4$ and
$R$.
By using
the last two observations, it holds that $u_{2,\lambda}\ge c_{10}>0$
on $(\lambda  D)\setminus D(\lambda,\delta_4,r_0)$. Therefore, by
\eqref{e:claim1} we get
 \bee\label{e:ubp}
\rho_{\lambda}(x)\ge u_{2, \lambda}(x) \ge \E_{x}\left[u_{2,
\lambda}\left(\wt X^{a,\lambda}_{\wt \tau^{a, \lambda}_{D(\lambda,
\delta_4,r_0)}} \right)\right] \ge c_{10} \P_{x}\left(\wt
X^{a,\lambda}_{\wt \tau^{a, \lambda}_{D( \lambda, \delta_4,r_0)}}
\in \lambda  D \right).
 \eee

Since the process $ \{\lambda (\wh X^{a}_{\lambda^{-2}t}-\wh X^{a}_{
0}), t\geq 0\}$ under $\P_x$ has the same distribution as $\{\wh
X^{a\lambda^{(\alpha-2)/\alpha}, \lambda}_t-\wh
X^{a\lambda^{(\alpha-2)/\alpha}, \lambda}_0, t\geq 0\}$ under
$\P_{\lambda x}$, we have from \eqref{e:lbp} that for $x \in D_Q(
\lambda^{-1} \delta_4, \lambda^{-1} r_0)$
\begin{eqnarray*}
&&\P_{x}\left(\wh X^{a}_{\tau^{a}_{ D_Q( \lambda^{-1} \delta_4,
\lambda^{-1} r_0 )}} \in  D \right)\\
&\ge&\P_{x}\left(\wh X^{a}_{\tau^{a}_{ D_Q( \lambda^{-1} \delta_4,
\lambda^{-1} r_0 )}} \in  D_Q(2\lambda^{-1} \delta_4, \lambda^{-1}
r_0 ) \setminus D_Q
( \lambda^{-1} \delta_4 ,\lambda^{-1} r_0 ) \right)\\
&=&\P_{\lambda x}\left(\wt X^{a,\lambda}_{\wt \tau^{a,
\lambda}_{ D(\lambda, \delta_4,r_0)}} \in D(\lambda,2\delta_4,r_0)
\setminus D(\lambda, \delta_4 ,r_0)
\right) \\
&\ge& c_9\rho_{\lambda}(\lambda x ) \,\ge\, c_{11} \delta_{
\lambda  D}(\lambda x ) \,=\,{c_{11}}{\lambda}
\delta_D (x)\, ,
\end{eqnarray*}
and, from \eqref{e:ubp}
\begin{eqnarray*}
\P_{x}\left(\wh X^{a}_{\tau^{a}_{ D_Q( \lambda^{-1} \delta_4,
\lambda^{-1} r_0 )}} \in  D \right) &=& \P_{\lambda x}\left(\wt
X^{a,\lambda}_{\wt \tau^{a, \lambda}_{
D(\lambda, \delta_4,r_0)}} \in \lambda  D \right)\\
 & \le& c_{12}\rho_{\lambda}(\lambda x ) \,\le\, c_{13}
\delta_{\lambda  D}(\lambda x ) \,=\,{c_{13}}{
\lambda} \delta_D (x).
\end{eqnarray*}
Finally by \eqref{e:scaling2_a} and \eqref{e:claim3},
\begin{eqnarray*}
\E_x\left[\wh \tau^{a}_{ D_Q( \lambda^{-1} \delta_4, \lambda^{-1}
r_0)}\right]&=&\int_{ D_Q( \lambda^{-1} \delta_4, \lambda^{-1} r_0 )} \wh
G^{a,1}_{ D_Q(
\lambda^{-1} \delta_4, \lambda^{-1} r_0 )} (x, y) dy\\
&=&\lambda^{d-2}\int_{ D_Q( \lambda^{-1} \delta_4, \lambda^{-1} r_0 )} \wh
G^{a\lambda^{(\alpha-2)/\alpha},\lambda}_{D(\lambda, \delta_4,r_0)}
(\lambda x, \lambda y)dy\\
&=&\lambda^{-2}\int_{D(\lambda, \delta_4,  r_0)} \wh
G^{a\lambda^{(\alpha-2)/\alpha},\lambda}_{D(\lambda, \delta_4,r_0)}
(\lambda x, z)dz\\
&=&\lambda^{-2} \E_{\lambda x} \left[\wt \tau^{a,
\lambda}_{D(\lambda, \delta_4,r_0)}
\right] \\
&\le &\lambda^{-2} \rho_{\lambda}(\lambda x ) \le  c_{14}
\lambda^{-2}  \delta_{\lambda D}(\lambda x )
\,=\,{c_{14}}{\lambda^{-1}} \delta_D (x).
\end{eqnarray*}
This completes the proof by taking $\delta_0=\delta_4$,
$C_8=c_{11}$, and $C_9=\max\{c_{13},c_{14}\}$. \qed

\medskip

We now derive exit distribution estimates for the process $X^a$
from those for $\wh X^a$ in Lemma \ref{L:2}.
Recall that $r_0=R/(4\sqrt{1+\Lambda^2})$.

\begin{lemma}\label{L:2.0_1}
There are  constants $\delta_0=\delta_0( R, M, \Lambda, \alpha)\in
(0, r_0)$, $C_8=C_8(R, M,\Lambda, \alpha)>0$ and $C_{10}=C_{10}(R,
M, \Lambda, \alpha)>0$ such that for every $a \in (0, M]$,
$\lambda\ge 1$,
$Q \in \partial D$ and  $x \in D_Q( \lambda^{-1} \delta_0 ,
\lambda^{-1} r_0 )$ with $\wt x =0$,
 \bee\label{e:L:2.0_1}
\P_{x}\left(X^{a}_{\tau^{a}_{   D_Q( \lambda^{-1} \delta_0 , \lambda^{-1}
r_0)}} \in D_Q(
 2\lambda^{-1} \delta_0 ,\lambda^{-1} r_0 )
\right) \ge\,{C_8}{\lambda} \delta_D (x),
 \eee
 \bee\label{e:L:2.0_2} \P_{x}\left(X^{a}_{\tau^{a}_{ D_Q(
\lambda^{-1} \delta_0, \lambda^{-1} r_0)}} \in  D\right)\le
{C_{10}}{\lambda} \delta_D (x)
 \eee
  and
 \bee\label{e:L:2.0_3} \E_x\left[\tau^a_{  D_Q(
\lambda^{-1} \delta_0, \lambda^{-1} r_0)}\right]\,\le\, {C_{10}}{\lambda^{-1}}
\delta_D (x).
 \eee
\end{lemma}

\pf Without loss of generality, we assume $Q=0$ and let $\phi:
\bR^{d-1}\to \bR$ be the $C^{1,1}$-function satisfying $\phi (\wt
0)= \nabla\phi (\wt 0)=0$, $\| \nabla \phi  \|_\infty \leq \Lambda$,
$| \nabla \phi (\wt y)-\nabla \phi (\wt z )| \leq \Lambda |\wt y-
\wt z|$ and $CS_Q$ be the corresponding coordinate system  such that
$$
B(Q, R)\cap D= \{(\wt y, \, y_d) \in CS_Q \cap B(0, R) : \ y_d >
\phi (\wt y) \}.
$$

Let $\delta_0$, $C_8$ and $C_9$ be the constants from the statement
of Lemma \ref{L:2}. Since $\mbox{diam}( D_Q( \lambda^{-1} \delta_0 ,
\lambda^{-1} r_0 )) \le \frac12$, we have that
$$
|x-y|^{-d-\alpha}\, \1_{\{|x-y|< 1\}} =|x-y|^{-d-\alpha} \quad
\text{for all } x, y\in D_Q(
 \lambda^{-1} \delta_0, \lambda^{-1} r_0).
$$

Let
$$
j(x) \,:= a^{\alpha} \sA(d, \alpha) |x|^{-(d+ \alpha)} \1_{\{|x| \ge
1\}}.
$$
Note that $\int_{\R^d} j(x)  dx <\infty$. Thus we can write $X^a_t
=\wh X^{a}_t +Z^a_t$ where $Z^a_t$ is a compound Poisson process
with the L\'evy density $j(x)$, independent of $\wh X^{a}_t$.
Since the jump size of $Z^a$ is greater than or equal to $1$ and
$\mbox{diam}( D_Q( \lambda^{-1} \delta_0, \lambda^{-1} r_0)) \le \frac12$, we
see from \eqref{e:L:3} that
$$
\E_x \left[\tau^{a}_{ D_Q( \lambda^{-1} \delta_0, \lambda^{-1} r_0)}\right]
\le\E_x \left[\wh \tau^{a}_{ D_Q( \lambda^{-1} \delta_0 , \lambda^{-1}
r_0)}\right]  \le {C_9}{\lambda^{-1}} \delta_D (x)\, .
$$
Moreover we have from \eqref{e:L:1} that
\begin{align*}
&\P_{x}\left(X^{a}_{\tau^{a}_{ D_Q( \lambda^{-1} \delta_0, \lambda^{-1} r_0)}}
\in D_Q( 2\lambda^{-1} \delta_0 ,\lambda^{-1} r_0 ) \right) \\=&\P_{x}\left(\wh
X^{a}_{\wh \tau^{a}_{ D_Q( \lambda^{-1} \delta_0, \lambda^{-1} r_0)}}
\in D_Q( 2\lambda^{-1} \delta_0 ,\lambda^{-1} r_0 ) \right) \ge
{C_8}{\lambda} \delta_D (x).
\end{align*}

We recall the notations from the proof of the previous lemma:
$$
\rho_\lambda(x):=y_d- \phi(\lambda^{-1} \wt y),
$$
$$
D(\lambda, r_1, r_2)\, :=\,\{ y\in CS_Q: r_1 >\rho_{\lambda}(y)
>0,\, |\wt y | < r_2 \},
$$
$$
\wt X^{a,\lambda} =\wh X^{a\lambda^{(\alpha-2)/\alpha},
\lambda} \qquad \hbox{and} \qquad  \wt \tau^{a,\lambda}_U:= \wh
\tau_U^{\, a\lambda^{(\alpha-2)/\alpha},\lambda}.
$$
Let  $a(\lambda):=a \lambda^{(\alpha-2)/\alpha}$, which is no larger
than $M$. By \eqref{e:levy},
\begin{align}
&\P_{x}\left(X^{a(\lambda)}_{\tau^{a(\lambda)}_{D(\lambda,
\delta_0,r_0)}} \in (\lambda  D)  \setminus D(\lambda, 2\delta_0
,2 r_0) \right)\label{e:ss22}\\
=&\ \E_x \left[ \int_0^{\tau^{a(\lambda)}_{D(\lambda,
\delta_0,r_0)}} \int_{ (\lambda  D)  \setminus D(\lambda,
2\delta_0 ,2 r_0)}\frac{a(\lambda)^{\alpha}\sA (d,
\alpha)}{|X^{a(\lambda)}_s-y|^{d+\alpha} }dy ds\right]
\nonumber\\
\le&\ c_1\sA (d, \alpha) (a(\lambda))^{\alpha} \left(\int_{(\lambda
D)  \setminus D(\lambda, 2\delta_0 ,2 r_0)} |y|^{-d-\alpha} dy \right)
\E_x \left[\tau^{a(\lambda)}_{D(\lambda, \delta_0,r_0)}\right]\nonumber\\
\le&\ c_2\sA (d, \alpha)  (a(\lambda))^{\alpha} \left(\int_{D(\lambda,
2\delta_0 ,r_0) \setminus D(\lambda, 3\delta_0/2
,r_0)}|y|^{-d-\alpha} dy \right)
\E_x \left[\tau^{a(\lambda)}_{D(\lambda, \delta_0,r_0)}\right]\nonumber\\
\le&\ c_3\E_x \left[ \int_0^{\tau^{a(\lambda)}_{D(\lambda,
\delta_0,r_0)}} \int_{ D(\lambda, 2 \delta_0,r_0) \setminus
D(\lambda, 3\delta_0/2 ,r_0)}\frac{a(\lambda)^{\alpha}\sA (d,
\alpha) }
{|X^{a(\lambda)}_s-y|^{d+\alpha} }dy ds\right]\nonumber\\
=&\ c_3\P_{x}\left(X^{a(\lambda)}_{\tau^{a(\lambda)}_{D(\lambda,
\delta_0, r_0)}} \in D(\lambda,2\delta_0,r_0) \setminus D( \lambda,
3\delta_0/2 ,r_0) \right).\nonumber
\end{align}

Thus by the above inequality and \eqref{e:ubp}, we have
\begin{align}
&\P_{x}\left(X^{a(\lambda)}_{ \tau^{a(\lambda)}_{D( \lambda,
\delta_0,r_0)}} \in \lambda  D \right)\label{e:ubp1}\\
=&\ \P_{x}\left(X^{a(\lambda)}_{ \tau^{a(\lambda)}_{D( \lambda,
\delta_0,r_0)}} \in (\lambda  D) \setminus D(\lambda, 2\delta_0
,2 r_0) \right)+\P_{x}\left(\wt X^{a,\lambda}_{\wt \tau^{a,
\lambda}_{D( \lambda, \delta_0,r_0)}} \in
D(\lambda, 2\delta_0 ,2 r_0) \right)\nonumber \\
\le&\ c_3   \P_{x}\left(X^{a(\lambda)}_{ \tau^{a(\lambda)}_{D(
\lambda, \delta_0,r_0)}} \in D(\lambda,2\delta_0,r_0) \setminus D(
\lambda, 3\delta_0/2 ,r_0) \right)+ \P_{x}\left(\wt
X^{a,\lambda}_{\wt \tau^{a, \lambda}_{D( \lambda, \delta_0,r_0)}}
\in \lambda  D \right)
\nonumber \\
= &\ c_3  \P_{x}\left(\wt X^{a,\lambda}_{\wt \tau^{a,
\lambda}_{D( \lambda, \delta_0,r_0)}}  \in D(\lambda,2\delta_0,r_0)
\setminus D( \lambda, 3\delta_0/2 ,r_0) \right)+ \P_{x}\left(\wt
X^{a,\lambda}_{\wt \tau^{a, \lambda}_{D( \lambda, \delta_0,r_0)}}
\in \lambda  D \right)\nonumber
 \\
\le &\ (c_3+1)  \P_{x}\left(\wt X^{a,\lambda}_{\wt \tau^{a,
\lambda}_{D( \lambda, \delta_0,r_0)}} \in \lambda  D \right) \le
c_4 \rho_\lambda(x). \nonumber
\end{align}
Since $(\lambda X^{a}_{\lambda^{-2} t}, t\geq 0)$ is the
independent sum of a Brownian motion and a symmetric $\alpha$-stable
process on $\R^d$ with weight $a(\lambda)$,
 we have from
\eqref{e:ubp1} that  for $x \in D_Q( \lambda^{-1} \delta_0, \lambda^{-1} r_0)$
\begin{align*}
\P_{x}\left(X^{a}_{\tau^{a}_{ D_Q( \lambda^{-1} \delta_0, \lambda^{-1} r_0)}}
\in  D \right) = &\ \P_{\lambda
x}\left(X^{a(\lambda)}_{\tau^{a(\lambda)}_{
D(\lambda, \delta_0,r_0)}} \in \lambda  D \right)\\
\le&\ c_4\rho_{\lambda}(\lambda x ) \,\le\, c_5
\delta_{\lambda  D}(\lambda x ) \,=\,{c_5}{ \lambda}
\delta_D (x).
\end{align*}
The proof is finished by taking $C_{10}=\max\{C_9, c_5\}$.
\qed

\section{Boundary Harnack principle}
In this section, we give the proof of the boundary Harnack principle for
the independent sum of a Brownian motion and a symmetric stable
process. We first prove the Carleson estimate for the independent
sum of a Brownian motion and a symmetric stable process on Lipschitz
open sets.

We recall that an open set $D$ in $\bR^d$ (when $d\ge 2$) is
said to be a Lipschitz open set if there exist a localization radius
$R_1>0$ and a constant $\Lambda_1 >0$ such that for every $Q\in
\partial D$, there exist a Lipschitz function $\phi=\phi_Q:
\bR^{d-1}\to \bR$ satisfying $\phi (0)= 0$, $| \phi (x)- \phi (y)|
\leq \Lambda_1 |x-y|$, and an orthonormal coordinate system $CS_Q$:
$y=(y_1, \dots, y_{d-1}, y_d)=:(\wt y, \, y_d)$ with its origin at
$Q$ such that
$$
B(Q, R_1)\cap D=\{ y=(\wt y, y_d)\in B(0, R_1) \mbox{ in } CS_Q: y_d
> \phi (\wt y) \}.
$$
The pair $(R_1, \Lambda_1)$ is called the characteristics of the Lipschitz
open set $D$.
Note that a
Lipschitz open set can be unbounded and disconnected.
For Lipschitz open set $D$ and every $Q\in \partial D$ and $ x \in
B(Q, R_1)\cap D$,
we define
$$
\rho_Q (x) := x_d -  \phi_Q (\tilde x)\, ,
$$
where $(\tilde x, x_d)$ is the coordinates of $x$ in $CS_Q$.

We recall that $X^a_t=X^0_t+a Y_t$ is a L\'evy process with
characteristic exponent $\Phi^a(x)=|x|^2+a^{\alpha}|x|^{\alpha}$.
This process may be obtained by subordinating a $d$-dimensional
Brownian motion
$W=(W_t,\, t\ge 0)$
by an independent subordinator
$T^a_t:=t+a^{2} T_t$ where $T=(T_t,\, t\ge 0)$ is an
 $\alpha/2$-stable subordinator. More precisely, the processes
$X^a_t$ and $W_{T^a_t}$ have the same distribution. Note that the
Laplace exponent corresponding to $T^a$ is equal to
$\phi^a(\lambda)=\lambda +a^{\alpha}\lambda^{\alpha/2}$. Let
$\sM_{\alpha/2}(t):=\sum_{n=0}^{\infty}(-1)^nt^{n\alpha/2}/\Gamma(1+n\alpha/2)$.
It follows by a straightforward integration that
$$
\int_0^{\infty}e^{-\lambda t}\sM_{1-\alpha/2}(a^{2\alpha/(2-\alpha)}
t)\, dt=\frac{1}{\phi^a(\lambda)}\, ,
$$
which shows that the potential density $u^a$ of the subordinator $T^a$
is given by
$$
u^a(t)=\sM_{1-\alpha/2}(a^{2\alpha/(2-\alpha)} t)\, .
$$
Since, for any $a>0$, $\phi^a$ is a complete Bernstein function, we
know that $u^a(\cdot)$ is a completely monotone function. In
particular, $u^a(\cdot)$ is a decreasing function. Since
$u^a(t)=u^1(a^{2\alpha/(2-\alpha)} t)$, we know that $a\mapsto
u^a(t)$ is a decreasing function.
 Therefore, if $0<a_1<a_2$, then $u^{a_1}(t)\ge u^{a_2}(t)$ for all
$t>0$. We will need this fact in the proof of next lemma.

\begin{lemma}\label{lower bound} Let $D\subset \R^d$ be a Lipschitz
open set with the characteristics $(R_1, \Lambda_1)$.
 There exists a constant $\delta=\delta(R_1, \Lambda_1, M)>0$ such
that for all $a\in [0,M]$ and $Q \in \partial D$, $x\in D$ with
$\rho_Q(x) < R_1/2$,
$$
\P_x(X^a_{\tau(x)}\in D^c)\ge \delta\, ,
$$
where $\tau(x):=\tau^{a}_{D\cap B(x,2\rho_Q(x))}=\inf\{t>0:\,
X^a_t\notin D\cap B(x,2\rho_Q(x))\}$.
\end{lemma}

\pf Clearly,
$$
\P_x\left(X^a_{\tau(x)}\in D^c\right)\,\ge\, \P_x\left(X^a_{\tau(x)}\in
D^c\cap B(x,2\rho_Q(x)) \right)\, \ge\, \P_x\left(X^a_{\tau(x)}\in \partial D\cap
B(x,2\rho_Q(x))\right).
$$
 Let $D_x:=D \cap B(x,2\rho_Q(x))$ and
$W^{D_x}$ be the subprocess of  Brownian motion $W$ killed upon leaving
 $D_x$. The process $Z^a$
defined by $Z^a_t:=W^{D_x}(T^a_t)$, where $T^a_t$ is an independent
subordinator described in the paragraph before the statement of the
lemma, is called a subordinate killed Brownian motion in $D_x$. We
will use $\zeta$ to denote the lifetime of $Z^a$.
 It is known from \cite{SV08} that
\begin{eqnarray*}
\P_x\left(X^a_{\tau(x)}\in \partial D\cap B(x,2\rho_Q(x))\right)&\ge&
\P_x\left(Z^a_{\zeta-}\in \partial D\cap B(x,2\rho_Q(x))\right)\\
&=&\E_x\left[u^a(\wt \tau_{D_x});\, W_{\wt \tau_{D_x}}\in
\partial D\cap B(x,2\rho_Q(x))\right]\, .
\end{eqnarray*}
Here and below, $\wt\tau_{U} :=\inf\{t>0:\, W_t\notin {U}\}$ is the exit
time of $W$ from ${U}$.
Denote $C_x:=\partial D\cap
B(x,2\rho_Q(x))$. Then
\begin{eqnarray}\label{ineq1}
\E_x\left[u^a(\wt\tau_{D_x}); \, W_{\wt\tau_{D_x}}\in C_x\right]
&\ge& \E_x\left[u^a(\wt\tau_{D_x}); \, W_{\wt\tau_{D_x}}\in C_x,
\wt\tau_{D_x} \le t\right]\\
&\ge & u^a(t)\P_x\left[W_{\wt\tau_{D_x}}\in C_x, \wt\tau_{D_x}
\le t\right]\nonumber \\
&\ge & u^a(t)\left(\P_x(W_{\wt\tau_{D_x}}\in
C_x)-\P_x(\wt\tau_{D_x}>t)\right)\, ,\nonumber
\end{eqnarray}
where $t>0$ will be chosen later.

Since $D$ is a Lipschitz open set with characteristics $(R_1,
\Lambda_1)$, there exist $\eta=\eta(\Lambda, R_1)>0$ and a cone
\begin{equation}\label{e:cone}
\sC:= \left\{y=(y_1, \dots ,y_d) \in \R^d :  y_d < 0 ,\, (y_1^2+ \cdots
+y_{d-1}^2 )^{1/2} < \eta |y_d| \right\}
\end{equation}
such that for every $z \in \partial D$, there is a cone $\sC_z$ with
vertex $z$, isometric to $\sC$, satisfying $\sC_z \cap B(Q,R_1)
\subset D^c$. Then by the scaling property of $W$ and symmetry
considerations, we have
\begin{eqnarray*}
\P_x(W_{\wt\tau_{D_x}}\in C_x) &\ge& \P_x\big( W_{\wt\tau_{B(x, 2
\rho_Q(x))}} \,\in\,
\partial B(x, 2 \rho_Q(x)) \cap \sC_{(\wt x,
\phi_Q(\wt x))} \big) \\
&\ge& \P_0\big( W_{\wt\tau_{B(0, 2) }} \,\in\, \partial B(0,2) \cap
( \sC + (\wt 0,  -1) ) \big),
\end{eqnarray*}
which is strictly positive.
Hence we can conclude that there exists $c_1=c_1(D)>0$ such that
\begin{equation}\label{ineq2}
\P_x(W_{\wt\tau_{D_x}}\in C_x)\ge c_1\, .
\end{equation}
Next,
\begin{equation}\label{ineq3}
\P_x(\wt\tau_{D_x} >t)\le \frac{\E_x[\wt\tau_{D_x}]}{t} \le
\frac{\E_x[\wt\tau_{B(x, 2\rho_Q(x))}]}{t} \le c_2
\frac{(\rho_Q(x))^2}{t}\le c_2\frac{R_1^2}{t}\, ,
\end{equation}
for some constant $c_2>0$. By using (\ref{ineq2}) and (\ref{ineq3})
in (\ref{ineq1}), we obtain that
$$
\E_x \left[u^a(\wt\tau_{D_x});\, W_{\wt\tau_{D_x}}\in C_x \right]
 \ge u^a(t)  \left(c_1-c_2\frac{R^2_1}{t}\right)\, .
$$
Now choose $t=t(R_1, \Lambda_1)>0$ large enough so that $c_1-c_2
R_1^2/t\ge c_1/2$. Then
$$
\E_x \left[u(\wt\tau_{D_x}); \, W_{\wt\tau_{D_x}}\in C_x\right]
\ge c_1 u^a(t)/2\ge c_1 u^M(t)/2=:\delta.
 $$
 The lemma is thus proved. \qed

Suppose that $D$ is an open set and that $U$ and $V$ are bounded open sets
with $V \subset \overline{V} \subset U$ and $ D \cap V \not= \emptyset$.
If $u$ vanishes continuously on $D^c\cap U$,
then by a finite covering argument, it is easy to see that $u$ is
bounded in an open neighborhood of $\partial D\cap V$.

\begin{lemma}\label{l:regularity}
Let $D$ be an open set and $U$ and $V$ be bounded open sets
with $V \subset \overline{V} \subset U$ and $ D \cap V \not= \emptyset$.
Suppose $u$ is a nonnegative function in $\R^d$ that is
harmonic in $D\cap U$ with respect to $X^a$ and vanishes
continuously on $D^c\cap U$. Then $u$ is regular harmonic in
$D\cap V$ with respect to $X^a$, i.e.,
\begin{equation}\label{e:regularity}
u(x)=\E_x\left[ u(X^a_{\tau^a_{D\cap V}})\right] \qquad \hbox{
for all }x\in D\cap V\, .
\end{equation}
\end{lemma}

\pf For $n\ge 1$, let  $B_n=\{y\in D\cap V:\, \delta_D(y)>1/n\} $.
Then for large $n$, $B_n$ is an non-empty open subset of $D\cap V$
whose closure is contained in $D\cap U$. Since $u$ is harmonic in
$D\cap U$ with respect to $X^a$, for $x\in D\cap V$ and $n$ large
enough so that $x\in B_n$,
 we have
that
$$
u(x)=\E_x\left[u\big(X^a_{\tau^a_{B_n}}\big)\right]=\E_x\left[
u\big(X^a_{\tau^a_{B_n}}\big); \, \tau^a_{B_n}< \tau^a_{D\cap V}\right]
+\E_x\left[u\big(X^a_{\tau^a_{B_n}}\big);\, \tau^a_{B_n}=\tau^a_{D\cap V}\right]\, .
$$
Hence
\begin{eqnarray}\label{e:reg2}
\lefteqn{\left|u(x)-\E_x\left[u\big(X^a_{\tau^a_{D \cap V}}\big)
\right] \right|}\\
&\le &\E_x\left[u\big(X^a_{\tau_{B_n}}\big);\, \tau^a_{B_n}<
\tau^a_{D\cap V}\right] +\E_x\left[u\big(X^a_{\tau^a_{D\cap
V}}\big);\,
\tau^a_{B_n}< \tau^a_{D\cap V}\right]\, .\nonumber
\end{eqnarray}
Since $\lim_{n \to  \infty} \tau^a_{B_n} =\tau^a_{D\cap V}$
$\P_x$-a.s., the second term in \eqref{e:reg2}
converges to $\E_x
\big[u\big(X^a_{\tau^a_{D\cap V}}\big);\, A \big]$ where
$A:=\cap_{n=1}^{\infty}\{\tau^a_{B_n} < \tau^a_{D\cap V}\}$.
 Note  that
$$X^a_{\tau^a_{D\cap V}}\in \partial
D\cap V \qquad \hbox{on } A.
$$
Hence  $u \big(X^a_{\tau^a_{D\cap V}}\big)=0$ on $A$, as
 $u$ is assumed to vanish on $D^c\cap
U$. Consequently
$$
\lim_{n \to \infty} \E_x\left[u\big(X^a_{\tau^a_{D\cap V}}\big);\,
\tau^a_{B_n}< \tau^a_{D\cap V}\right]=0\, .
$$
For the first term in \eqref{e:reg2}, note that
$\delta_D(X^a_{\tau_{B_n}}) \le 1/n$ on $\{\tau^a_{B_n}<
\tau^a_{D\cap V}\}$. Therefore, by the assumption that $u$
vanishes continuously on $D^c\cap U$, one has $\lim_{n \to
\infty} u(X^a_{\tau_{B_n}})=0$. Moreover, since $u$ vanishes
continuously on $(\partial D)\cap U$, there is $n_0\geq 1$ so
that $u$ is bounded in $D\cap V\setminus B_{n_0}$. So by the
bounded convergence theorem we have
$$
\lim_{n \to  \infty}
\E_x\left[u\big(X^a_{\tau_{B_n}}\big);\, \tau^a_{B_n}< \tau^a_{D\cap
V}\right]=0\,.
$$
This proves the lemma. \qed

\noindent {\bf Proof of Proposition \ref{uhp}}.
We know from the parabolic Harnack inequality from Theorem 6.7 of \cite{CK08}
that Harnack inequality holds for the process $X:=X^1$.
That is,
there exists a constant $c_1=c_1(\alpha, M)>0$
such that for any $r\in (0, M^{\alpha/(2-\alpha)}]$, $x_0\in \R^d$ and
any function $v \ge 0$ harmonic in
$B(x_0, r)$ with respect to
 $X$,
 we have
\begin{equation}\label{e:HP1}
v(x)\le c_1v(y) \qquad \mbox{ for all } x, y\in B(x_0, \frac{r}2).
\end{equation}
Now the proposition is an easy consequence of \eqref{e:HP1}. In
fact, note that  for any $a\in (0, M]$, $X^a$ has the same
distribution as $\lambda X_{\lambda^{-2}t}$, where $\lambda =
a^{\alpha/(\alpha-2)} \geq M^{\alpha/(\alpha-2)}$. Consequently, if
$u$ is harmonic in $B(x_0, r)$ with respect to  $X^a$ where $r \in
(0,1]$, then $v(x):=u(\lambda x)$ is harmonic in $B(\lambda^{-1}x_0,
\lambda^{-1}r)$ with respect to $X$ and $\lambda^{-1}r \le
M^{\alpha/(2-\alpha)}$. So by \eqref{e:HP1}
$$
u(\lambda x)=v(x)\le c_1v(y)= c_1u(\lambda y)
\qquad \mbox{ for all } x, y\in B(\lambda^{-1}x_0, \lambda^{-1} r/2).
$$
That is,
$$
u( x)\le c_1u( y) \qquad \mbox{ for all } x, y\in B(x_0, r/2).
$$
\qed

\begin{thm}[Carleson estimate]\label{carleson}
Let $D\subset \R^d$ be a Lipschitz open set with the characteristics
$(R_1, \Lambda_1)$. Then there exists a positive
 constant $A=A(\alpha, \Lambda_1, R_1, M)$ such that for $a \in (0,
M]$, $Q\in \partial D$, $0<r<R_1/2$, and any nonnegative function
$u$ in $\R^d$ that is harmonic in $D \cap B(Q, r)$ with respect to
$X^{a}$ and vanishes continuously on $ D^c \cap B(Q, r)$, we have
\begin{equation}\label{e:carleson}
u(x)\le A u(x_0) \qquad \hbox{for }  x\in D\cap B(Q,r/2).
\end{equation}
where $x_0\in D
\cap B(Q,r)$ with $\rho_Q(x_0)=r/2$.
\end{thm}

\pf Fix $a\in (0, M]$. Since $D$ is Lipschitz and $r<R_1/2$, by the
uniform Harnack principle in
 Proposition
\ref{uhp} and a standard chain
argument, it suffices to prove (\ref{e:carleson}) for $x\in D\cap
B(Q,r/12)$ and $\wt x_0 = \wt Q$. Without loss of generality, we may
assume that $u(x_0)=1$.
In this proof, constants $\delta, \beta, \eta$ and $c_i$'s are always independent of $r$ and $a$.

Choose $0<\gamma <\alpha/(d+\alpha)$ and let
$$
B_0=D\cap B(x,2\rho_Q(x))\, ,\qquad B_1=B(x,r^{1-\gamma}
\rho_Q(x)^{\gamma})\, .
$$
Further, set
$$
B_2=B(x_0,\rho_Q(x_0)/3)\, ,\qquad B_3=B(x_0, 2\rho_Q(x_0)/3)
$$
and
$$
 \tau_0=\inf\{t>0:\, X^a_t\notin B_0\}\, , \qquad  \tau_2=\inf\{t>0:\,
X^a_t\notin B_2\}.
 $$
 By Lemma \ref{lower bound}, there exists
$\delta=\delta(R_1, \Lambda_1, M)>0$ such that
\begin{equation}\label{e:c:1}
\P_x(X^a_{\tau_0}\in D^c)\ge \delta\, ,\quad x\in B(Q,r/4)\, .
\end{equation}
By the uniform Harnack principle in
  Proposition
  \ref{uhp} and a chain
argument, there exists $\beta$ such that
\begin{equation}\label{e:c:2}
u(x)<(\rho_Q(x)/r)^{-\beta} u(x_0)\, ,\quad x\in D\cap B(Q,r/4)\, .
\end{equation}
In view of Lemma \ref{l:regularity}, $u$ is regular harmonic in $B_0$
with respect to $X^a$. So
\begin{equation}\label{e:c:3}
u(x)=\E_x\big[u\big(X^a_{\tau_0}\big); X^a_{\tau_0}\in B_1\big]+
\E_x\big[u\big(X^a_{\tau_0}\big); X^a_{\tau_0}\notin B_1\big]
\qquad \hbox{for } x\in B(Q, r/4) .
\end{equation}
We first show that there exists $\eta>0$ such that
\begin{equation}\label{e:c:4}
\E_x\big[u\big(X^a_{\tau_0}\big); X^a_{\tau_0}\notin B_1\big]\le u(x_0)  \quad
\hbox{if } x\in D \cap B(Q,r/12) \hbox{ with } \rho_Q(x) < \eta r\,
.
\end{equation}
Let $\eta_0 :=2^{-2(d+\alpha)/d }$, then for $\rho_Q(x)< \eta_0r$,
$$ (\rho_Q(x))^{d /(\alpha+d)} < 1/4 \quad \hbox{ and } \quad
 2\rho_Q(x) \le r^{1-\gamma} \rho_Q(x)^{\gamma} - 2\rho_Q(x).
$$
  Thus if $x\in D
\cap B(Q,r/12)$ with $\rho_Q(x) < \eta_0r$,  then $|x-y|\le 2|z-y|$
for $z\in B_0$, $y\notin B_1$. Thus we have by \eqref{e:levy} and
Lemma  \ref{L:2.00}
\begin{eqnarray}\label{e:c:5}
&&\E_x\big[u\big(X^a_{\tau_0}\big); X^a_{\tau_0}\notin B_1\big]
\\
&=&\sA (d, \alpha) \int_{B_0}G^a_{B_0}(x,z)\int_{|y-x|>r^{1-\gamma}
\rho_Q(x)^{\gamma}}a^{\alpha}|z-y|^{-d-\alpha}u(y)\, dy\, dz
\nonumber \\
&\le & 2^{d+\alpha}\sA (d, \alpha) \int_{B_0}G^a_{B_0}(x,z)dz
\int_{|y-x|>r^{1-\gamma}\rho_Q(x)^{\gamma}}
a^{\alpha}|x-y|^{-d-\alpha}u(y)\,
dy \nonumber \\
&\le&  2^{d+\alpha}\sA (d, \alpha)  \E_x[\tau_{B(x, 2\rho_Q(x))}]
\int_{|y-x|>r^{1-\gamma}\rho_Q(x)^{\gamma}}
a^{\alpha}|x-y|^{-d-\alpha}u(y)\, dy\nonumber \\
&\le &2^{d+\alpha}\sA (d, \alpha)  c_1 \rho_Q(x)^2
\left(\int_{|y-x|>r^{1-\gamma}\rho_Q(x)^{\gamma}, |y-x_0|>2
\rho_Q(x_0)/3}
a^{\alpha}|x-y|^{-d-\alpha}u(y)\, dy \right.\nonumber \\
&  &\quad \quad \quad \quad \quad +\left.\int_{|y-x_0|\le
2\rho_Q(x_0)/3}a^{\alpha}
|x-y|^{-d-\alpha}u(y)\, dy\right)\nonumber \\
&=:& c_2 \rho_Q(x)^2 (I_1+I_2)\, .\nonumber
\end{eqnarray}
On the other hand, for $z\in B_2$ and $y\notin
B_3$, we have $|z-y|\le |z-x_0|+|x_0-y|\le \rho_Q(x_0)/3+|x_0-y|\le
2|x_0-y|$. we have again by \eqref{e:levy} and Lemma
\ref{L:2.00}
\begin{eqnarray}\label{e:c:6}
u(x_0)&\ge& \E_{x_0}\left[u(X^a_{\tau_2}), X^a_{\tau_2}\notin B_3\right]
\\
&\ge & \sA (d, \alpha)
\int_{B_2}G^a_{B_2}(x_0,z)\int_{|y-x_0|>2\rho_Q(x_0)/3}
a^{\alpha}|z-y|^{-d-\alpha}u(y)\, dy \, dz\nonumber \\
&\ge &2^{-d-\alpha}\sA (d, \alpha) \int_{B_2}G^a_{B_2}
(x_0,z)dz\int_{|y-x_0|>2\rho_Q(x_0)/3}
a^{\alpha}|x_0-y|^{-d-\alpha}u(y)\, dy \nonumber \\
&\ge &2^{-d-\alpha}\sA (d, \alpha)  c_3 (\rho_Q(x_0)/3)^2
\int_{|y-x_0|>2\rho_Q(x_0)/3}
a^{\alpha}|x_0-y|^{-d-\alpha}u(y)\, dy \nonumber \\
&=&c_4 \rho_Q(x_0)^2 \int_{|y-x_0|>2\rho_Q(x_0)/3}
a^{\alpha}|x_0-y|^{-d-\alpha}u(y)\, dy\, .\nonumber
\end{eqnarray}

Suppose now that $|y-x|\ge r^{1-\gamma}\rho_Q(x)^{\gamma}$ and $x\in
B(Q,r/4)$. Then
$$|y-x_0|\le |y-x|+r\le
|y-x|+r^{\gamma}\rho_Q(x)^{-\gamma}|y-x|\le
2r^{\gamma}\rho_Q(x)^{-\gamma}|y-x|.
$$
Therefore
\begin{eqnarray}\label{e:c:7}
\lefteqn{I_1=\int_{|y-x|>r^{1-\gamma}\rho_Q(x)^{\gamma},
|y-x_0|>2\rho_Q(x_0)/3}a^{\alpha} |x-y|^{-d-\alpha}u(y)\, dy}
\\
&\le & \int_{|y-x_0|>2\rho_Q(x_0)/3}
(2^{-1}(\rho_Q(x)/r)^{\gamma})^{-d-\alpha}
a^{\alpha}|y-x_0|^{-d-\alpha} u(y)\, dy\nonumber \\
&= &
2^{d+\alpha}(\rho_Q(x)/r)^{-\gamma(d+\alpha)}\int_{|y-x_0|>2\rho_Q(x_0)/3}
a^{\alpha}|y-x_0|^{-d-\alpha}u(y)\, dy\nonumber \\
&\le & 2^{d+\alpha}(\rho_Q(x)/r)^{-\gamma(d+\alpha)} c_4^{-1}
\rho_Q(x_0)^{-2} u(x_0)\nonumber \\
&=& c_5 (\rho_Q(x)/r)^{-\gamma(d+\alpha)}\rho_Q(x_0)^{-2}u(x_0)\,
,\nonumber
\end{eqnarray}
 where the last inequality is due to \eqref{e:c:6}.

If $|y-x_0|<2\rho_Q(x_0)/3$, then $|y-x|\ge |x_0-Q|-|x-Q|-|y-x_0|
>\rho_Q(x_0)/6$.
This together with the uniform Harnack principle in
 Proposition
\ref{uhp}
implies that
\begin{eqnarray}\label{e:c:8}
\lefteqn{I_2 =\int_{|y-x_0|\le 2\rho_Q(x_0)/3}a^{\alpha}
|x-y|^{-d-\alpha} u(y)\, dy}\\
&\le &c_6 \int_{|y-x_0|\le 2\rho_Q(x_0)/3}a^{\alpha}
|x-y|^{-d-\alpha} u(x_0)\, dy\nonumber \\
&\le & c_6 u(x_0)\int_{|y-x|>\rho_Q(x_0)/6} a^{\alpha}
|x-y|^{-d-\alpha}\, dy \,=\, c_7 a^{\alpha}\rho_Q(x_0)^{-\alpha}
u(x_0)\, .\nonumber
\end{eqnarray}

Combining (\ref{e:c:5})-(\ref{e:c:8}) we obtain
\begin{eqnarray}\label{e:c:9}
&&\E_x[u(X^a_{\tau_0});\, X^a_{\tau_0}\notin B_1]\\
&\le & c_2\rho_Q(x)^2\left(c_5 (\rho_Q(x)/r)^{-\gamma( d +
\alpha)}\rho_Q(x_0)^{-2}u(x_0)+ c_7 a^{\alpha}
\rho_Q(x_0)^{-\alpha} u(x_0)\right)\nonumber \\
&\le &c_8 u(x_0)\left(\rho_Q(x)^2(\rho_Q(x)/
r)^{-\gamma(d+\alpha)}\rho_Q(x_0)^{-2}
+a^{\alpha}\rho_Q(x)^2\rho_Q(x_0)^{-\alpha}\right)\nonumber \\
&\le & c_9 u(x_0) \left((\rho_Q(x)/r)^{2- \gamma(d + \alpha)}
+a^{\alpha}\rho_Q(x)^2 r^{-\alpha}\right)\, ,\nonumber
\end{eqnarray}
where in the last inequality we used the fact that $\rho_Q(x_0)=r/2$.
Choose now
$\eta\in (0, \eta_0)$ so that
$$
c_9\,\left(\eta^{2-\gamma(d+\alpha)}+\eta^2 M^{\alpha}\right)\,\le\,
1\, .
$$
Then  for $x\in  D \cap
B(Q,r/12)$ with $\rho_Q(x) < \eta r$, we have by \eqref{e:c:9}
\begin{eqnarray*}
\E_x\left[u(X^a_{\tau_0});\, X^a_{\tau_0}\notin B_1\right]\,&\le\,&
c_9\, u(x_0)\left(\eta^{2-\gamma(d+\alpha)}+\eta^2
r^{2-\alpha}M^{\alpha}\right) \\
&\le\, &c_9 \left(\eta^{2-\gamma(d+\alpha)} +
\eta^2 M^{\alpha}\right) u(x_0)\le u(x_0)\, .
\end{eqnarray*}

We now prove the Carleson estimate \eqref{e:carleson} for $x\in
D\cap B(Q, r/12)$ by a method of contradiction. Recall that
 $u(x_0)=1$. Suppose that there exists $x_1\in D\cap B(x,r/12)$ such
that $u(x_1)\ge K>\eta^{-\beta}\vee (1+\delta^{-1})$, where $K$ is a
 constant to be specified later. By \eqref{e:c:2} and the
 assumption $u(x_1)\ge K>\eta^{-\beta}$, we have
$(\rho_Q(x_1)/r)^{-\beta}>u(x_1)\ge K> \eta^{-\beta}$, and hence
$\rho_Q(x_1)<\eta r$. Let $B_0$, $B_1$ and $\tau_0$ be now defined
with respect to the point $x_1$ instead of $x$. Then by
(\ref{e:c:3}),  (\ref{e:c:4}) and $K>1+\delta^{-1}$,
$$
K\le u(x_1)\le \E_{x_1}\left[u(X^a_{\tau_0}); X^a_{\tau_0} \in
B_1\right]+1\, ,
$$
and hence
$$
\E_{x_1}\left[u(X^a_{\tau_0}); X^a_{\tau_0}\in B_1\right] \ge
u(x_1)-1 > \frac{1}{1+\delta}\, u(x_1)\, .
$$
In the last inequality of the display above we used the assumption that
 $u(x_1)\ge K>1+\delta^{-1}$. If $K \ge 2^{\beta/\gamma}$, then $D^c\cap
B_1\subset D^c \cap B(Q,r)$. By using the assumption that $u=0$ on
$D^c\cap B(Q, r)$, we get from \eqref{e:c:1}
$$
\E_{x_1}[u(X^a_{\tau_0}), X^a_{\tau_0}\in B_1]=\E_{x_1}[
u(X^a_{\tau_0}), X^a_{\tau_0}\in B_1\cap D] \le \P_x(X^a_{\tau_0}\in
D) \, \sup_{B_1}u  \le (1-\delta) \, \sup_{B_1}u \, .
$$
Therefore, $\sup_{B_1} u> u(x_1)/((1+\delta)(1-\delta))$, i.e.,
there exists a point $x_2\in D$ such that
$$
|x_1-x_2|\le r^{1-\gamma}\rho_Q(x_1)^{\gamma} \quad \hbox{ and }
\quad
u(x_2)>\frac{1}{1-\delta^2}\, u(x_1)\ge \frac{1}{1-\delta^2}\, K\, .
$$
By induction, if $x_k\in D\cap B(Q, r/12)$ with $u(x_k)\geq
K/(1-\delta^2)^{k-1}$ for $k\ge 2$, then there exists $x_{k+1}\in D$
such that
\begin{equation}\label{e:c:10}
|x_k-x_{k+1}|\le r^{1-\gamma}\rho_Q(x_k)^{\gamma}  \quad \hbox{ and }
\quad
u(x_{k+1}) > \frac{1}{1-\delta^2}\, u(x_k)>
\frac{1}{(1-\delta^2)^k}\, K\, .
\end{equation}
{}From (\ref{e:c:2}) and (\ref{e:c:10}) it follows that
$\rho_Q(x_{k})/r \le (1-\delta^2)^{(k-1)/\beta}K^{-1/\beta}$, for
every $k\ge 1$. Therefore,
\begin{eqnarray*}
|x_k-Q|&\le&|x_1-Q|
+\sum_{j=1}^{k-1}|x_{j+1}-x_j|\le \frac{r}{12} +
 \sum_{j=1}^{\infty} r^{1-\gamma}\rho_Q(x_j)^{\gamma}\\
&\le &\frac{r}{12}+r^{1-\gamma}\sum_{j=1}^{\infty}(1-\delta^2)^{(j-1)
\gamma/\beta}K^{-\gamma/\beta}r^{\gamma}\\
&=&\frac{r}{12}+r^{1-\gamma}r^{\gamma}K^{-\gamma/\beta}
\sum_{j=0}^{\infty}(1-\delta^2)^{j\gamma/\beta}\\
&=&\frac{r}{12}+ r K^{-\gamma/\beta}\,
\frac{1}{1-(1-\delta^2)^{\gamma/\beta}}.
\end{eqnarray*}
Choose $K=\eta\vee (1+\delta^{-1})\vee
12^{\beta/\gamma}(1-(1-\delta^2)^{\gamma/\beta})^{-\beta/\gamma}$.
Then $K^{-\gamma/\beta}\, (1-(1-\delta^2)^{\gamma/\beta})^{-1}\le
1/12$, and hence $x_k\in D\cap B(Q,r/6)$ for every $k\ge 1$. Since
$\lim_{k\to \infty}u(x_k)=+\infty$, this contradicts the fact that
$u$ is bounded on $B(Q,r/2)$. This contradiction shows that $u(x)<
K$ for every $x\in D\cap B(Q, r/12)$. This completes the proof of
the theorem. \qed

\noindent {\bf Proof of Theorem \ref{t:main} }. Since $D$ is a
$C^{1,1}$ open set and $r<R$,
by the uniform Harnack principle in
  Proposition
 \ref{uhp} and a standard
chain argument, it suffices to prove \eqref{e:bhp_m} for $x,y \in D
\cap B(Q,rr_0/8)$.
In this proof, constants $\eta$ and $c_i$'s are always independent of $r$ and $a$.

We recall  that
$r_0=\frac{R}{4\sqrt{1+\Lambda^2}}$
and $\delta_0\in
(0, r_0)$ is the constant in the statement of Lemma \ref{L:2.0_1}.

For any $r\in (0, R]$ and
 $x\in D\cap B(Q, rr_0/8)$, let $Q_x$ be
the point $Q_x \in \partial D$ so that $|x-Q_x|=\delta_{D}(x)$ and
let $x_0:=Q_x+\frac{r}{8}(x-Q_x)/|x-Q_x|$. We choose a
$C^{1,1}$-function $\phi: \bR^{d-1}\to \bR$ satisfying $\phi
(0)= \nabla\phi (0)=0$, $\| \nabla \phi  \|_\infty \leq \Lambda$, $|
\nabla \phi (y)-\nabla \phi (z)| \leq \Lambda |y-z|$, and an
orthonormal coordinate system $CS$ with its origin at $Q_x$ such
that
$$
B(Q_x, R)\cap D=\{ y=(\wt y, y_d) \in B(0, R) \mbox{ in } CS: y_d >
\phi (\wt y) \}.
$$
In the coordinate system $CS$ we have $\wt x = \wt 0$ and $x_0=(\wt
0, r/8)$. For any $b_1, b_2>0$, we define
$$
D(b_1, b_2):=\left\{ y=(\wt y, y_d) \mbox{ in } CS: 0<y_d-\phi(\wt
y)<b_1r\delta_0/8, \ |\wt y| < b_2rr_0/8 \right\}.
$$
It is easy to see that
$D(2, 2)\subset D\cap B(Q, r/2)$.
In fact, since $r_0 \le \frac{1}{8\Lambda}$ and  $\delta_0 \le
\frac{1}{8\Lambda}$,
 for every $z \in  D(2, 2)$
$$
|z-Q| \le |Q-x|+|x-Q_x|+|Q_x-z| \le
 \frac{r}8+\frac{r}8
+ |z_d- \phi(\wt
z)|+ |\phi(\wt z)| < \frac{r}{4}
 ( 1+ \delta_0) \le \frac{r}{2}.
$$
Thus if $u$ is a nonnegative function on $\R^d$ that is harmonic in
$D\cap B(Q, r)$ with respect to $X^a$ and vanishes continuously in
$D^c\cap B(Q, r)$, then, by Lemma \ref{l:regularity}, $u$ is regular
harmonic in $D\cap B(Q,r/2)$ with respect to $X^a$, hence also in
$D(2, 2)$.
Thus by the uniform Harnack principle in
 Proposition
   \ref{uhp}, we have
\begin{eqnarray}
u(x) &= & \E_x\left[u\big(X^a_{\tau^a_{ D(1,1)}}\big)\right] \ge
\E_x\left[u\big(X^a_{\tau^a_{ D(1,1)}}\big); X^a_{
\tau^a_{ D(1,1)}} \in  D(2,1)\right]\label{e:BHP2}\\
&\ge& c_1 u(x_0) \P_x\Big( X^a_{\tau^a_{ D(1,1)}} \in D (2,1)\Big)
\ge c_2 u(x_0) \delta_D(x) /r.\nonumber
\end{eqnarray}
In the last inequality above we have used \eqref{e:L:2.0_1}.

 Let $w=(\wt 0, rr_0/16)$.
  Then it is easy to see that there
exists a constant $\eta=\eta(\Lambda, r_0, \delta_0)\in (0, 1)$ such
that $B(w, \eta rr_0/16)\in D(1, 1)$.  By \eqref{e:levy} and
Lemma \ref{L:2.00},
\begin{eqnarray*}
u(w) &\ge& \E_{w}\left[u\big(X^a_{\tau^a_{ D(1,1)}}\big);
X^a_{\tau^a_{ D(1,1)}} \notin  D(2,2)\right]\\
&=& \sA (d, \alpha) a^{\alpha}\int_{ D(1,1)} G^a_{ D(1,1)}(w,z)
\int_{\R^d\setminus D(2,2)} \frac{u(y)}{|z-y|^{d+\alpha}}dydz\\
&\ge& c_3 a^{\alpha}\E_{w}\big[\tau^a_{B(w, \eta rr_0/(16)
)}\big] \int_{\R^d\setminus  D(2,2)}
\frac{u(y)}{|w-y|^{d+\alpha}}dy\\
&\ge& c_4 a^{\alpha} r^2   \int_{\R^d\setminus  D(2,2)}
\frac{u(y)}{|w-y|^{d+\alpha}}dy.
\end{eqnarray*}

Hence by \eqref{e:L:2.0_3},
\begin{eqnarray*}
&&\E_{x}\left[u \left(X^a_{\tau^a_{ D(1,1)}}\right); \,
X^a_{\tau^a_{ D(1,1)}}
\notin  D(2,2)\right]\\
&=& \sA (d, \alpha) a^{\alpha} \int_{ D(1,1)} G^a_{ D(1,1)}(x,z)
\int_{\R^d\setminus
D(2,2)}  \frac{u(y)}{|z-y|^{d+\alpha}}dydz\\
&\le& c_5 a^{\alpha}\E_x[\tau^a_{ D(1,1)}]  \int_{\R^d\setminus
D(2,2)}  \frac{u(y)}{|w-y|^{d+\alpha}}dy\\
&\le& c_6 a^{\alpha} \delta_D(x) r  \int_{\R^d\setminus  D(2,2)}
\frac{u(y)}{|w-y|^{d+\alpha}}dy \,\leq \,  \frac{c_6 \,
\delta_D(x)}{c_4 \, r} u(w).
\end{eqnarray*}

On the other hand, by the uniform Harnack principle
 (Proposition
  \ref{uhp}) and the Carleson estimate (Theorem \ref{carleson}), we
have
$$
\E_x\left[u\left(X^a_{\tau^a_{ D(1,1)}}\right);\,  X^a_{\tau^a_{
D(1,1)}} \in D(2,2)\right] \,\le\, c_7 \, u(x_0) \P_x\left(
X^a_{\tau^a_{D (1,1)}} \in  D(2,2)\right)\,\le\, c_8 \, u(x_0)
\delta_D(x) /r.
$$
In the last inequality above we have used \eqref{e:L:2.0_2}.
Combining the two inequalities above, we get
\begin{eqnarray}
  u(x)
 &=& \, \E_x\left[u \left(X^a_{\tau^a_{ D(1,1)}}\right);
\, X^a_{\tau^a_{ D(1,1)}} \in D(2,2)\right]
 \label{e:BHP1} \\
&& +\E_x\left[ u \left(X^a_{\tau^a_{ D(1,1)}}\right); \,
X^a_{\tau^a_{ D(1,1)}} \notin  D(2,2)\right] \nonumber \\
&\le& \, \frac{c_8}{r} \delta_D(x) u(x_0)   + \frac{c_6 \,
\delta_D(x)}{c_4\, r} u(w) \nonumber \\
&\le&  \, \frac{c_{9}}{r}\, \delta_D(x) (u(x_0)  + u(w))\nonumber \\
& \le & \frac{c_{10}}{r}\, \delta_D(x) u(x_0)  .\nonumber
\end{eqnarray}
In the last inequality above we have used the
uniform Harnack principle
 (Proposition \ref{uhp}).

From \eqref{e:BHP2}-\eqref{e:BHP1}, we have that for every $x, y\in
  D \cap B(Q, rr_0/8)$,
$$
\frac{u(x)}{u(y)}\,\le \,
\frac{c_{10}}{c_2}\,\frac{\delta_D(x)}{\delta_D(y)},
$$
which  proves the theorem. \qed

\vskip 0.3truein

{\bf Zhen-Qing Chen}

Department of Mathematics, University of Washington, Seattle,
WA 98195, USA

E-mail: \texttt{zchen@math.washington.edu}

\bigskip

{\bf Panki Kim}

Department of Mathematical Sciences and Research Institute of Mathematics,
Seoul National University,
San56-1 Shinrim-dong Kwanak-gu,
Seoul 151-747, Republic of Korea

E-mail: \texttt{pkim@snu.ac.kr}

\bigskip

{\bf Renming Song}

Department of Mathematics, University of Illinois, Urbana, IL 61801, USA

E-mail: \texttt{rsong@math.uiuc.edu}

\bigskip

{\bf
Zoran Vondra\v{c}ek}

Department of Mathematics,
University of Zagreb,
Bijeni\v{c}ka c.~30,
Zagreb, Croatia

Email: \texttt{vondra@math.hr}

\end{document}